\documentclass[a4paper,11pt,twoside,reqno]{amsart}

\usepackage[utf8]{inputenc}
\usepackage[plainpages=false,pdfpagelabels=true]{hyperref}
\usepackage{amssymb,amsthm,slashed}
\usepackage[margin=1in]{geometry}
\usepackage[all]{xy}
\usepackage{dsfont}

\newtheorem{Satz}{Theorem}[section]
\newtheorem{Prop}[Satz]{Proposition}
\newtheorem{Lem}[Satz]{Lemma}

\newtheorem{Cor}[Satz]{Corollary}

\theoremstyle{definition}
\newtheorem{Dfn}[Satz]{Definition}
\newtheorem{Bem}[Satz]{Remark}

\newcommand{\sff}{\mathrm{I\!I}}

\allowdisplaybreaks[1]

\renewcommand{\epsilon}{\varepsilon}

\newcommand{\R}{\ensuremath{\mathbb{R}}}

\numberwithin{equation}{section}

\title{Energy methods for Dirac-type equations in two-dimensional Minkowski space}
\author{Volker Branding}
\date{\today}
\address{University of Vienna, Faculty of Mathematics\\
Oskar-Morgenstern-Platz 1, 1090 Vienna, Austria\\}
\email{volker.branding@univie.ac.at}
\begin{document}

\begin{abstract}
In this article we develop energy methods 
for a large class of linear and nonlinear Dirac-type 
equations in two-dimensional Minkowski space.
We will derive existence results for several Dirac-type
equations originating in quantum field theory, in particular
for Dirac-wave maps to compact Riemannian manifolds.
\end{abstract} 

\maketitle

\section{Introduction and results}
In quantum field theory spinors are used to describe fermions,
which are elementary particles of half-integer spin.
The equations that govern their behavior are both linear and nonlinear Dirac equations. 
Linear Dirac equations are employed to model free fermions.
However, to model the interaction of fermions one has to take into
account nonlinearities.

In mathematical terms spinors are sections in the spinor bundle, which is
a vector bundle defined on the underlying manifold. 
Its existence requires the vanishing of the second Stiefel-Whitney class,
which is a topological condition.
The natural operator acting on spinors is the Dirac operator,
which is a first-order differential operator. If the underlying 
manifold is Riemannian the Dirac operator is elliptic, if the manifold 
is Lorentzian the Dirac operator is hyperbolic.
The fact that the Dirac operator is of first order usually leads to
technical problems since there are less tools available compared to
second order operators such as the Laplacian.

In the case of a Riemannian manifold many results on the qualitative behavior
of nonlinear Dirac equations have been obtained recently, see
\cite{MR2390834,MR3558358,branding2016nonlinear}. However, obtaining an existence result 
for nonlinear Dirac equations in the Riemannian case is rather complicated.

This article is supposed to be the first step to develop energy methods for linear and nonlinear
Dirac equations on Lorentzian manifolds in a geometric framework.
For an introduction to linear geometric wave equations on globally hyperbolic manifolds we refer to
\cite{MR2298021}.

As a starting point we will stick to the case of
two-dimensional Minkowski space, the generalization to higher dimensional globally
hyperbolic manifolds will be treated in a subsequent work.
Most of the analytic results on nonlinear Dirac equation in Minkowski space make use of a global trivialization
of the spinor bundle and investigate the resulting coupled system of partial differential equations
of complex-valued functions. In our approach we do not make use of a global trivialization of the spinor bundle,
but derive energy estimates for the spinor itself. This approach
seems to be the natural one from a geometric point of view.

This article is organized as follows: After presenting the necessary background on spin geometry
in two-dimensional Minkowski space in the next subsection, we will focus on the analysis 
of linear Dirac equations in section 2. 
Afterwards, in section 3 we will consider several models from quantum field theory that involve
nonlinear Dirac equations. Making use of the energy methods developed before we 
derive existence results for some of these models.
The last section is devoted to the study of Dirac-wave maps from two-dimensional Minkowski space
taking values in a compact Riemannian manifold. Again, by application of suitable energy methods,
we are able to derive an existence result for the latter.

\subsection{Spin geometry in two-dimensional Minkowski space}
First, let us fix the notations that will be used throughout this article.
In addition, we want to recall several facts on spin geometry in the setting of a Lorentzian manifold.
In this article, we will make use of the Einstein summation convention, that is we sum over repeated indices.

In the following we will consider two-dimensional Minkowski space \(\R^{1,1}\) with metric \((+,-)\) 
and global coordinates \((t,x)\). The tangent vectors of \(\R^{1,1}\) will be denoted by \(\partial_t,\partial_x\).

The spinor bundle over \(\R^{1,1}\) will be denoted by \(\Sigma\R^{1,1}\) and sections in this bundle will be called 
\emph{spinors}. Note that \(\Sigma\R^{1,1}\) can be globally trivialized, however we will not often make use of this fact.
On \(\Sigma\R^{1,1}\) there exists a metric connection and we have a Hermitian, but indefinite, scalar product
denoted by \(\langle\cdot,\cdot\rangle\).
We will use the convention that this scalar product is linear in the first and anti-linear in the second slot.

We denote Clifford multiplication of a spinor \(\psi\) with a tangent vector \(X\) by \(X\cdot\psi\).
Note that in contrast to the Riemannian case, Clifford multiplication is symmetric, that is
\begin{align*}
\langle X\cdot\xi,\psi\rangle_{\Sigma\R^{1,1}}=\langle\xi,X\cdot\psi\rangle_{\Sigma\R^{1,1}}
\end{align*}
for all \(\psi,\xi\in\Gamma(\Sigma\R^{1,1})\) and all \(X\in T\R^{1,1}\).
In addition, the Clifford relations
\begin{align*}
X\cdot Y\cdot\psi+Y\cdot X\cdot\psi=-2g(X,Y)\psi
\end{align*}
hold for all \(X,Y\in T\R^{1,1}\), where \(g\) is the metric of \(\R^{1,1}\).

The Dirac operator on two-dimensional Minkowski is defined as (with \(\epsilon_j=g(e_j,e_j)=\pm 1\))
\begin{align*}
D:=\sum_{j=1}^2\epsilon_je_j\cdot\nabla_{e_j}=\partial_t\cdot\nabla_{\partial_t}-\partial_x\cdot\nabla_{\partial_x}.
\end{align*}
Here, \(\nabla\) denotes the connection on \(\Sigma\R^{1,1}\) and 
\(e_i,i=1,2\) is a pseudo-orthonormal basis of \(T\R^{1,1}\).

Note that, in contrast to the Riemannian case, the Dirac operator defined above is not formally 
self-adjoint with respect to the \(L^2\)-norm but satisfies
\[
\int_{\R^{1,1}}\langle\xi,D\psi\rangle d\mu=-\int_{\R^{1,1}}\langle D\xi,\psi\rangle d\mu
\]
for all \(\psi,\xi\in\Sigma\R^{1,1}\).
For this reason, we will mostly consider the operator \(iD\), since this combination is formally self-adjoint with
respect to the \(L^2\)-norm.

For many of the analytic questions discussed in this article it will be necessary to have 
a positive-definite scalar product on the spinor bundle in order to establish energy estimates.

For this reason we consider the positive definite scalar product
\begin{align*}
\langle\partial_t\cdot,\rangle,
\end{align*}
where \(\partial_t\) denotes the globally-defined timelike vector field.
The resulting norm will be denoted by \(||_\beta\), that is
\begin{align*}
0\leq|\psi|^2_\beta:=\langle\partial_t\cdot\psi,\psi\rangle
\end{align*}
for \(\psi\in\Gamma(\Sigma\R^{1,1})\).

For more details on spin geometry on Lorentzian manifolds we refer to \cite{MR701244}
and also \cite{MR3074852,MR2121740,MR3551253}.

\begin{Bem}
Note that we have two kinds of natural scalar products on the spinor bundle
in the semi-Riemannian case. On the one hand we have the (geometric) scalar product
that is invariant under the Spin group, but indefinite. On the other hand we have 
the (analytic) scalar product, which is positive definite but breaks the geometric invariance.
\end{Bem}

\section{Linear Dirac equations in two-dimensional Minkowski space}
In this section we derive conserved energies for solutions of linear Dirac-type equations.
Later on, we will generalize these methods to the non-linear case.

We start be analyzing solutions of
\begin{align}
\label{psi-kernel}
D\psi=0.
\end{align}
Note that, for a solution of \eqref{psi-kernel}, we have the following identities
\begin{align}
\label{identity-tt-xx-psi-kernel}\frac{\partial}{\partial t}\langle\partial_t\cdot\psi,\psi\rangle-\frac{\partial}{\partial x}\langle\partial_x\cdot\psi,\psi\rangle
=&\langle D\psi,\psi\rangle+\langle\psi, D\psi\rangle=0, \\
\label{identity-tx-xt-psi-kernel}\frac{\partial}{\partial t}\langle\partial_x\cdot\psi,\psi\rangle-\frac{\partial}{\partial x}\langle\partial_t\cdot\psi,\psi\rangle
=&\langle(\partial_t\cdot\nabla_{\partial_x}-\partial_x\cdot\nabla_{\partial_t})\psi,\psi\rangle
-\langle\psi,(\partial_t\cdot\nabla_{\partial_x}+\partial_x\cdot\nabla_{\partial_t})\psi\rangle=0.
\end{align}

We will use these identities to derive several conservation laws.

\begin{Lem}
\label{lemma1-energies-psi-kernel}
Let \(\psi\in\Gamma(\Sigma\R^{1,1})\) be a solution of \eqref{psi-kernel}.
Then the energies 
\begin{align*}
E_1(t)=&\frac{1}{2}\int_\R|\psi|^2_\beta dx,\\
E_2(t)=&\frac{1}{2}\int_\R\big(\big|\frac{\partial}{\partial x}|\psi|^2_\beta\big|^2+\big|\frac{\partial}{\partial t}|\psi|^2_\beta\big|^2\big)dx, \\
E_3(t)=&\frac{1}{2}\int_\R(|\nabla_{\partial_t}\psi|_\beta^2+|\nabla_{\partial_x}\psi|_\beta^2)dx
\end{align*}
are conserved.
\end{Lem}
\begin{proof}
The first claim follows from integrating \eqref{identity-tt-xx-psi-kernel}.
Differentiating \eqref{identity-tt-xx-psi-kernel} with respect to \(t\) and \eqref{identity-tx-xt-psi-kernel} with respect to \(x\) we find
\begin{align*}
\frac{\partial^2}{\partial t^2}\langle\partial_t\cdot\psi,\psi\rangle=\frac{\partial^2}{\partial x\partial t}\langle\partial_x\cdot\psi,\psi\rangle
=\frac{\partial^2}{\partial x^2}\langle\partial_t\cdot\psi,\psi\rangle.
\end{align*}
Hence, \(|\psi|^2_\beta\) solves the one-dimensional wave equation, which yields the second statement.
The third assertion follows since \(\psi\) solves \(\nabla^2_{\partial_t}\psi=\nabla^2_{\partial_x}\psi\).
\end{proof}

We can use \eqref{identity-tt-xx-psi-kernel} and \eqref{identity-tx-xt-psi-kernel} to find conserved energies that involve higher \(L^p\) norms
of \(\psi\).

\begin{Lem}
\label{lemma2-energies-psi-kernel}
Let \(\psi\in\Gamma(\Sigma\R^{1,1})\) be a solution of \eqref{psi-kernel}. Then the energy
\begin{align*}
E_4(t)=\int_\R(|\psi|^4_\beta+|\langle\partial_x\cdot\psi,\psi\rangle|^2)dx
\end{align*}
is conserved.
\end{Lem}
\begin{proof}
Making use of \eqref{identity-tt-xx-psi-kernel} and \eqref{identity-tx-xt-psi-kernel} we calculate
\begin{align*}
\frac{d}{dt}\frac{1}{2}\int_\R|\psi|^4_\beta dx=&\int_\R|\psi|_\beta^2\frac{\partial}{\partial t}\langle\partial_t\cdot\psi,\psi\rangle dx\\
=&\int_\R|\psi|_\beta^2\frac{\partial}{\partial x}\langle\partial_x\cdot\psi,\psi\rangle dx \\
=&-\int_\R\langle\partial_x\cdot\psi,\psi\rangle\frac{\partial}{\partial x}\langle\partial_t\cdot\psi,\psi\rangle dx \\
=&-\int_\R\langle\partial_x\cdot\psi,\psi\rangle\frac{\partial}{\partial t}\langle\partial_x\cdot\psi,\psi\rangle dx \\
=&-\frac{d}{dt}\frac{1}{2}\int_\R|\langle\partial_x\cdot\psi,\psi\rangle|^2_\beta dx,
\end{align*}
which proves the claim.
\end{proof}

\begin{Bem}
It is straightforward to generalize the previous conservation law to any \(L^p\) norm of \(\psi\) making use 
of \eqref{identity-tt-xx-psi-kernel} and \eqref{identity-tx-xt-psi-kernel}.
\end{Bem}

\begin{Prop}
Let \(\psi\in\Gamma(\Sigma\R^{1,1})\) be a solution of \eqref{psi-kernel}.
Then the following energy 
\begin{align*}
E_5(t)=\frac{1}{2}\int_\R\big(\big|\frac{\partial}{\partial x}e(\psi)\big|^2+\big|\frac{\partial}{\partial t}e(\psi)\big|^2\big)dx
\end{align*}
is conserved, where
\[
e(\psi):=\frac{1}{2}(|\nabla_{\partial_t}\psi|_\beta^2+|\nabla_{\partial_x}\psi|_\beta^2).
\]
\end{Prop}
\begin{proof}
By a direct calculation we find that \(e(\psi)\) solves the one-dimensional wave equation,
which yields the statement.
\end{proof}

\begin{Bem}
By the Sobolev embedding \(H^1\hookrightarrow L^\infty\) this yields a pointwise bound on \(e(\psi)\).
It is straightforward to also bound higher derivatives of \(\psi\).
\end{Bem}

As a next step we investigate if the same conservation laws still hold if \(\psi\) solves 
a linear Dirac equation with a right hand side.
To this end, we consider the linear Dirac equation
\begin{align}
\label{psi-linear}
iD\psi=\lambda\psi,\qquad \lambda\in\R.
\end{align}
Note that \eqref{psi-linear} arises as critical point of the functional
\begin{align*}
S(\psi)=\int_{\R^{1,1}}(\langle\psi,iD\psi\rangle-\lambda|\psi|^2)d\mu,
\end{align*}
which leads to the prefactor of \(i\) in front of the Dirac operator.
Moreover, note that we do not use the definite \(|\psi|^2_\beta\)-norm.

For a solution of \eqref{psi-linear}, we have the following identities
\begin{align}
\label{identity-tt-xx-psi-linear}\frac{\partial}{\partial t}\langle\partial_t\cdot\psi,\psi\rangle-\frac{\partial}{\partial x}\langle\partial_x\cdot\psi,\psi\rangle
=&-\lambda(\langle i\psi,\psi\rangle+\langle\psi, i\psi\rangle=0, \\
\label{identity-tx-xt-psi-linear}\frac{\partial}{\partial t}\langle\partial_x\cdot\psi,\psi\rangle-\frac{\partial}{\partial x}\langle\partial_t\cdot\psi,\psi\rangle
=&\langle(\partial_t\cdot\nabla_{\partial_x}-\partial_x\cdot\nabla_{\partial_t})\psi,\psi\rangle
-\langle\psi,(\partial_t\cdot\nabla_{\partial_x}+\partial_x\cdot\nabla_{\partial_t})\psi\rangle \\
=\nonumber&2\lambda\langle i\partial_x\cdot\partial_t\cdot\psi,\psi\rangle.
\end{align}

\begin{Lem}
Let \(\psi\in\Gamma(\Sigma\R^{1,1})\) be a solution of \eqref{psi-linear}.
Then the energies 
\begin{align*}
E_1(t)=&\frac{1}{2}\int_\R|\psi|^2_\beta dx,\\
E_6(t)=&\frac{1}{2}\int_\R(|\nabla_{\partial_t}\psi|_\beta^2+|\nabla_{\partial_x}\psi|_\beta^2-\lambda^2|\psi|^2_\beta)dx
\end{align*}
are conserved.
\end{Lem}
\begin{proof}
The first statement follows from integrating \eqref{identity-tt-xx-psi-linear}.
Regarding the second claim, we note that due to the prefactor of \(i\) we obtain the following wave-type equation
when squaring the Dirac operator
\[
\nabla^2_{\partial_t}\psi-\nabla^2_{\partial_x}\psi=\lambda^2\psi,
\]
which yields the second statement.
\end{proof}

\begin{Prop}
For given initial data \(\psi(0,x)=\psi_0(x)\in H^1(\R,\Sigma\R^{1,1})\) the solution of \eqref{psi-linear} exists globally in that space.
\end{Prop}
\begin{proof}
Making use of the conserved energy \(E_6(t)\) we obtain the following energy inequality
\begin{align*}
\int_\R|\partial_x|\psi|_\beta|^2dx\leq \int_\R|\nabla_{\partial_x}\psi|_\beta^2dx\leq \lambda^2\int_\R|\psi|^2_\beta dx+E_6(t)\leq C
\end{align*}
for a uniform constant \(C\).
Moreover, by the Sobolev embedding \(H^1\hookrightarrow L^\infty\) this yields a pointwise bound
on \(|\psi|_\beta\). Consequently, the solution of \eqref{psi-linear} exists globally.
\end{proof}

\subsection{Twisted spinors}
\label{subsection-twisted-spinors}
In this subsection we want to discuss if the previous results still hold when we consider
twisted spinors, which are sections in the spinor bundle that is twisted by some additional vector bundle \(F\).

To this end let \(F\) be a Hermitian vector bundle with a metric connection.
Moreover, we will assume that we have a positive definite scalar product on \(F\).
On the twisted bundle \(\Sigma\R^{1,1}\otimes F\) we obtain a metric connection 
induced from the connections on \(\Sigma\R^{1,1}\) and \(F\), which we will denote by \(\tilde{\nabla}\),
via setting
\begin{align*}
\tilde{\nabla}:=\nabla^{\Sigma\R^{1,1}}\otimes\mathds{1}^{F}+\mathds{1}^{\Sigma\R^{1,1}}\otimes\nabla^{F}.
\end{align*}

The twisted Dirac operator \(D^F\colon\Gamma(\Sigma\R^{1,1}\otimes F)\to\Gamma(\Sigma\R^{1,1}\otimes F)\) is defined by
\begin{align*}
D^F:=\epsilon_je_j\cdot\tilde{\nabla}_{e_j}=\partial_t\cdot\tilde{\nabla}_{\partial_t}-\partial_x\cdot\tilde{\nabla}_{\partial_x}.
\end{align*}
In contrast to the spinor bundle \(\Sigma\R^{1,1}\) over Minkowski space the vector bundle \(F\) is not supposed
to be flat such that it may have non-vanishing curvature. We will denote its curvature endomorphism by \(R^F(\cdot,\cdot)\).

\begin{Lem}
The square of the twisted Dirac operator \(D^F\) satisfies the following Weitzenboeck formula
\begin{align}
\label{weitzenboeck-twisted}
(D^F)^2=-\tilde{\nabla}^2_{\partial_t}+\tilde{\nabla}^2_{\partial_x}-\partial_t\cdot\partial_x\cdot R^F(\partial_t,\partial_x),
\end{align}
where \(R^F\) denotes the curvature of the vector bundle \(F\).
\end{Lem}
\begin{proof}
We calculate
\begin{align*}
(D^F)^2=&(\partial_t\cdot\tilde{\nabla}_{\partial_t}-\partial_x\cdot\tilde{\nabla}_{\partial_x})(\partial_t\cdot\tilde{\nabla}_{\partial_t}-\partial_x\cdot\tilde{\nabla}_{\partial_x}) \\
=&-\tilde{\nabla}^2_{\partial_t}+\tilde{\nabla}^2_{\partial_x}-\partial_t\cdot\partial_x\cdot\tilde{\nabla}_{\partial_t}\tilde{\nabla}_{\partial_x}	
-\partial_x\cdot\partial_t\cdot\tilde{\nabla}_{\partial_x}\tilde{\nabla}_{\partial_t}\\
=&-\tilde{\nabla}^2_{\partial_t}+\tilde{\nabla}^2_{\partial_x}-\partial_t\cdot\partial_x\cdot R^F(\partial_t,\partial_x),
\end{align*}
which completes the proof.
\end{proof}
Note that, compared to the Riemannian case, we have a different sign in front of the curvature term of the vector bundle \(F\).
In addition, we do not get a scalar curvature contribution in  \eqref{weitzenboeck-twisted} since we are restricting
ourselves to two-dimensional Minkowski space.

\begin{Bem}
Most of the Dirac type equations studied in quantum field theory involve twisted Dirac operators \cite{MR677998}.
In particular, the spinors that are considered in the standard model of elementary particle physics
are sections in the spinor bundle twisted by some vector bundle.
\end{Bem}

Again, we start by deriving several energy estimates for solutions of 
\begin{equation}
\label{psi-kernel-twisted}
D^F\psi=0.
\end{equation}

For solutions of \eqref{psi-kernel-twisted} we obtain the following identities
\begin{align}
\label{identity-tt-xx-psi-kernel-twisted} \frac{\partial}{\partial t}\langle\partial_t\cdot\psi,\psi\rangle-\frac{\partial}{\partial x}\langle\partial_x\cdot\psi,\psi\rangle
=&\langle D^F\psi,\psi\rangle+\langle\psi, D^F\psi\rangle=0, \\
\label{identity-tx-xt-psi-kernel-twisted}\frac{\partial}{\partial t}\langle\partial_x\cdot\psi,\psi\rangle-\frac{\partial}{\partial x}\langle\partial_t\cdot\psi,\psi\rangle
=&\langle(\partial_t\cdot\tilde{\nabla}_{\partial_x}-\partial_x\cdot\tilde{\nabla}_{\partial_t})\psi,\psi\rangle
-\langle\psi,(\partial_t\cdot\tilde{\nabla}_{\partial_x}+\partial_x\cdot\tilde{\nabla}_{\partial_t})\psi\rangle=0.
\end{align}

\begin{Prop}
\label{prop-psi-kernel-twisted-a}
Let \(\psi\in\Gamma(\Sigma\R^{1,1}\otimes F)\) be a solution of \eqref{psi-kernel-twisted}.
Then the energies 
\begin{align*}
\tilde{E}_1(t)=&\frac{1}{2}\int_R|\psi|^2_\beta dx,\\
\tilde{E}_2(t)=&\frac{1}{2}\int_\R\big(\big|\frac{\partial}{\partial x}|\psi|^2_\beta\big|^2+\big|\frac{\partial}{\partial t}|\psi|^2_\beta\big|^2\big)dx
\end{align*}
are conserved.
\end{Prop}
\begin{proof}
This follows as in the proof of Lemma \ref{lemma1-energies-psi-kernel}.
\end{proof}

\begin{Lem}
\label{lemma2-energies-psi-kernel-twisted}
Let \(\psi\in\Gamma(\Sigma\R^{1,1}\otimes F)\) be a solution of \eqref{psi-kernel-twisted}. Then the energy
\begin{align*}
\tilde{E}_4(t)=\int_\R(|\psi|^4_\beta+|\langle\partial_x\cdot\psi,\psi\rangle|^2)dx
\end{align*}
is conserved.
\end{Lem}
\begin{proof}
This follows as in the proof of Lemma \ref{lemma2-energies-psi-kernel} making use of \eqref{identity-tt-xx-psi-kernel-twisted} and 
\eqref{identity-tx-xt-psi-kernel-twisted}.
\end{proof}

\begin{Bem}
Again, it is straightforward to also find conserved energies involving higher \(L^p\) norms of \(\psi\)
for solutions of \eqref{psi-kernel-twisted}.
\end{Bem}

We set 
\begin{align*}
\tilde{E}_3(t):=\frac{1}{2}\int_\R(|\tilde{\nabla}_{\partial_t}\psi|_\beta^2+|\tilde{\nabla}_{\partial_x}\psi|_\beta^2)dx.
\end{align*}

When we try to control derivatives of solutions of \eqref{psi-kernel-twisted} it will be
necessary to control the curvature of the vector bundle \(F\).

\begin{Lem}
Let \(\psi\in\Gamma(\Sigma\R^{1,1}\otimes F)\) be a solution of \eqref{psi-kernel-twisted}.
Then \(\tilde{E}_3(t)\) satisfies 
\begin{align*}
\frac{d}{dt}\tilde{E}_3(t)\leq\tilde{E}_3(t)+\frac{\big||\psi|^2_\beta\big|_{L^\infty}}{2}\int_\R|R^F(\partial_t,\partial_x)|^2 dx.
\end{align*}

\end{Lem}
\begin{proof}
By assumption we have \(D^F\psi=0\) and consequently also \((D^F)^2\psi=0\). Now, we calculate
\begin{align*}
\frac{d}{dt}\tilde{E}_3(t)
=&\int_\R(\langle\partial_t\cdot\tilde{\nabla}^2_{\partial_t}\psi,\tilde{\nabla}_{\partial_t}\psi\rangle 
+\langle\partial_t\cdot\tilde{\nabla}_{\partial_t}\tilde{\nabla}_{\partial_x}\psi,\tilde{\nabla}_{\partial_x}\psi\rangle) dx\\
=&\int_\R(\langle\partial_t\cdot(\tilde{\nabla}^2_{\partial_t}-\tilde{\nabla}^2_{\partial_x})\psi,\tilde{\nabla}_{\partial_t}\psi\rangle 
+\langle\partial_t\cdot R^F(\partial_t,\partial_x)\psi,\tilde{\nabla}_{\partial_x}\psi\rangle) dx \\
=&\int_\R(\langle\partial_x\cdot R^F(\partial_t,\partial_x)\psi,\tilde{\nabla}_{\partial_t}\psi\rangle
+\langle\partial_t\cdot R^F(\partial_t,\partial_x)\psi,\tilde{\nabla}_{\partial_x}\psi\rangle) dx \\
\leq& \tilde{E}_3(t)+\frac{\big||\psi|^2_\beta\big|_{L^\infty}}{2}\int_\R|R^F(\partial_t,\partial_x)|^2dx
\end{align*}
yielding the result.
\end{proof}

As a next step we discuss if the previous methods can still be employed if \(\psi\) solves a linear Dirac equation
with non-trivial right-hand side, that is
\begin{align}
\label{psi-twisted-linear}
iD^F\psi=\lambda\psi,\qquad \lambda\in\R.
\end{align}

For a solution of \eqref{psi-twisted-linear}, we have the following identities
\begin{align}
\label{identity-tt-xx-psi-twisted-linear}\frac{\partial}{\partial t}\langle\partial_t\cdot\psi,\psi\rangle-\frac{\partial}{\partial x}\langle\partial_x\cdot\psi,\psi\rangle
=&0, \\
\label{identity-tx-xt-psi-twisted-linear}\frac{\partial}{\partial t}\langle\partial_x\cdot\psi,\psi\rangle-\frac{\partial}{\partial x}\langle\partial_t\cdot\psi,\psi\rangle
=&2\lambda\langle i\partial_x\cdot\partial_t\cdot\psi,\psi\rangle.
\end{align}

\begin{Lem}
Let \(\psi\in\Gamma(\Sigma\R^{1,1}\otimes F)\) be a solution of \eqref{psi-twisted-linear}.
Then the energy
\begin{align*}
\tilde{E}_1(t)=&\frac{1}{2}\int_\R|\psi|^2_\beta dx
\end{align*}
is conserved.
\end{Lem}
\begin{proof}
This follows by integrating \eqref{identity-tt-xx-psi-twisted-linear}.
\end{proof}

\begin{Lem}
Let \(\psi\in\Gamma(\Sigma\R^{1,1}\otimes F)\) be a solution of \eqref{psi-twisted-linear}.
Then the following inequality holds
\begin{align*}
\int_\R(|\psi|^4_\beta+|\langle\partial_x\cdot\psi,\psi\rangle|^2)dx\leq Ce^{|\lambda|t},
\end{align*}
where the positive constant \(C\) depends on \(\psi_0\).
\end{Lem}
\begin{proof}
The proof follows by a direct calculation making use of \eqref{identity-tt-xx-psi-twisted-linear} and \eqref{identity-tx-xt-psi-twisted-linear}.
\end{proof}

\begin{Lem}
Let \(\psi\in\Gamma(\Sigma\R^{1,1}\otimes F)\) be a solution of \eqref{psi-twisted-linear}.
Then the following inequality holds
\begin{align*}
\frac{d}{dt}\int_\R(|\tilde{\nabla}_{\partial_t}\psi|_\beta^2+|\tilde{\nabla}_{\partial_x}\psi|_\beta^2-\lambda^2|\psi|^2_\beta)dx
\leq\int_\R(|\tilde{\nabla}_{\partial_t}\psi|_\beta^2+|\tilde{\nabla}_{\partial_x}\psi|_\beta^2)dx
+\int_\R|\psi|^2_\beta|R^F(\partial_t,\partial_x)|^2dx.
\end{align*}
\end{Lem}
\begin{proof}
By the Weitzenboeck formula \eqref{weitzenboeck-twisted} we find
\begin{align*}
\tilde{\nabla}^2_{\partial_t}\psi-\tilde{\nabla}^2_{\partial_x}\psi=\lambda^2\psi-\partial_t\cdot\partial_x\cdot R^F(\partial_t,\partial_x)\psi.
\end{align*}
In addition, we calculate
\begin{align*}
\frac{d}{dt}\frac{1}{2}\int_\R(|\tilde{\nabla}_{\partial_t}\psi|_\beta^2+&|\tilde{\nabla}_{\partial_x}\psi|_\beta^2-\lambda^2|\psi|^2_\beta)dx \\
=&\int_\R(\langle\partial_x\cdot R^F(\partial_t,\partial_x)\psi,\tilde{\nabla}_{\partial_t}\psi\rangle
+\langle\partial_t\cdot R^F(\partial_t,\partial_x)\psi,\tilde{\nabla}_{\partial_x}\psi\rangle) dx
\end{align*}
yielding the result.
\end{proof}

\begin{Bem}
We have to impose a bound of the form \(\int_\R|\psi|_\beta^2|R^F(\partial_t,\partial_x)|^2dx\leq C\)
if we want to deduce an energy estimate for solutions of \eqref{psi-twisted-linear}.
\end{Bem}

\section{Nonlinear Dirac equations in two-dimensional Minkowski space}
In this section we want to investigate if the energy methods developed for linear Dirac equations
in the previous section can also be applied to the nonlinear case. 
The equations we will study mostly arise in quantum field theory, however, 
some of them also are connected to problems in differential geometry, see for example \cite{MR2424048}.

Up to now there exist many analytic results on nonlinear Dirac equations in two-dimensional
Minkowski space. A general framework for semilinear hyperbolic systems was developed in \cite{MR2257730},
for a recent survey on nonlinear Dirac equations see \cite{MR3383628} and references therein.

\subsection{The Thirring Model}
First, we will focus on a famous model from quantum field theory, the \emph{Thirring model}.
This model was introduced in \cite{THIRRING195891} to describe the self-interaction of a Dirac field in two-dimensional
Minkowski space.
In the physics literature there exists a huge number of results on the Thirring model
and also in the mathematical literature many results, including existence results, have been established.
Most of the methods employed so far make use of a global trivialization of the spinor bundle
over two-dimensional Minkowski space yielding existence results for the Thirring model, 
see for example \cite{MR2802088,MR2330348,MR932441,MR3325607,MR0463658}.
Making use of our energy methods we will also provide an existence result for the Thirring model.

The action functional for the Thirring model is given by
\begin{align}
\label{functional-thirring-geometric}
S(\psi)=\int_{\R^{1,1}}(\langle\psi,iD\psi\rangle-\lambda|\psi|^2-\frac{\kappa}{2}\epsilon_j\langle\psi,e_j\cdot\psi\rangle\langle\psi,e_j\cdot\psi\rangle)d\mu
\end{align}
with real parameters \(\kappa\) and \(\lambda\). In physics, \(\lambda\) is usually interpreted as mass,
whereas \(\kappa\) describes the strength of interaction.
The critical points of \eqref{functional-thirring-geometric} are given by
\begin{align}
\label{psi-thirring-geometric}
iD\psi=\lambda\psi+\kappa\epsilon_j\langle\psi,e_j\cdot\psi\rangle e_j\cdot\psi.
\end{align}

\begin{Bem}
Note that for \(\lambda=0\) solutions of \eqref{psi-thirring-geometric}
are invariant under scaling, that is if \(\psi\) is a solution of \eqref{psi-thirring-geometric},
then 
\begin{align*}
\psi(t,x)\to r\psi(r^2t,r^2x),
\end{align*}
where \(r\) is a positive number, is also a solution.
\end{Bem}

For a solution of \eqref{psi-thirring-geometric} the following identity holds
\begin{align*}
\partial_x\cdot\nabla_{\partial_t}\psi-\partial_t\cdot\nabla_{\partial_x}
=&i\lambda\partial_x\cdot\partial_t\cdot\psi+i\epsilon_j\kappa\langle\psi,e_j\cdot\psi\rangle\partial_x\cdot\partial_t\cdot e_j\cdot\psi\\
=&i\lambda\partial_x\cdot\partial_t\cdot\psi-i\kappa\langle\psi,\partial_t\cdot\psi\rangle\partial_x\cdot\psi
+i\kappa\langle\psi,\partial_x\cdot\psi\rangle\partial_t\cdot\psi
\end{align*}
leading to the two equations
\begin{align}
\label{identity-tt-xx-psi-thirring-geometric}
\frac{\partial}{\partial t}\langle\partial_t\cdot\psi,\psi\rangle-\frac{\partial}{\partial x}\langle\partial_x\cdot\psi,\psi\rangle
=&\lambda(\langle i\psi,\psi\rangle+\langle\psi,i\psi\rangle)
-\kappa\langle\psi,e_j\cdot\psi\rangle(\langle ie_j\cdot\psi,\psi\rangle+\langle\psi,ie_j\cdot\psi\rangle)=0, \\
\label{identity-tx-xt-psi-thirring-geometric}
\frac{\partial}{\partial t}\langle\partial_x\cdot\psi,\psi\rangle-\frac{\partial}{\partial x}\langle\partial_t\cdot\psi,\psi\rangle=&2\lambda\langle i\partial_x\cdot\partial_t\cdot\psi,\psi\rangle. 
\end{align}
Hence, for a solution of \eqref{psi-thirring-geometric} the energy
\begin{align}
E_1(t)=&\frac{1}{2}\int_\R|\psi|^2_\beta dx
\end{align}
is conserved again. 
By combining \eqref{identity-tt-xx-psi-thirring-geometric} and \eqref{identity-tx-xt-psi-thirring-geometric} we find
\begin{align*}
\Box |\psi|^2_\beta=2\lambda\frac{\partial}{\partial x}\langle i\partial_x\cdot\partial_t\cdot\psi,\psi\rangle 
=4\lambda(\langle i\nabla_{\partial_t}\psi,\psi\rangle+2|\psi|^2_\beta(\lambda-\kappa|\psi|^2)).
\end{align*}
This identity turns out to be very useful in the case of the massless Thirring model,
that is for solutions of \eqref{psi-thirring-geometric} with \(\lambda=0\).
More precisely, we find
\begin{Prop}
Let \(\psi\in\Gamma(\Sigma\R^{1,1})\) be a solution of \eqref{psi-thirring-geometric} with \(\lambda=0\). 
Then
\begin{align}
|\psi|^2_\beta \leq C
\end{align}
for a positive constant \(C\).
\end{Prop}
\begin{proof}
Since \(|\psi|^2_\beta\) solves the one-dimensional wave equation, we again get a conserved energy 
and the result follows from the Sobolev embedding \(H^1\hookrightarrow L^\infty\).
\end{proof}

In order to treat the massive Thirring model we need several auxiliary lemmata.
\begin{Lem}
Let \(\psi\in\Gamma(\Sigma\R^{1,1})\) be a solution of \eqref{psi-thirring-geometric}. Then the following wave type equation holds
\begin{align}
\label{thirring-wave}
\nabla^2_{\partial_t}\psi-\nabla^2_{\partial_x}\psi=&\lambda^2\psi+\epsilon_j\lambda\kappa\langle\psi,e_j\cdot\psi\rangle e_j\cdot\psi 
+2\kappa\lambda\langle i\partial_x\cdot\partial_t\cdot\psi,\psi\rangle i\partial_x\cdot\partial_t\cdot\psi \\
&\nonumber+\kappa\langle\psi,\partial_t\cdot\psi\rangle(-2i\nabla_{\partial_t}\psi-\lambda\partial_t\cdot\psi)
+\kappa\langle\psi,\partial_x\cdot\psi\rangle(-2i\partial_t\cdot\partial_x\cdot\nabla_{\partial_t}\psi-\lambda\partial_x\cdot\psi) \\
&\nonumber+\kappa^2(|\langle\psi,\partial_t\cdot\psi\rangle|^2\psi+|\langle\psi,\partial_x\cdot\psi\rangle|^2\psi
-2\langle\psi,\partial_x\cdot\psi\rangle\langle\psi,\partial_t\cdot\psi\rangle\partial_x\cdot\partial_t\cdot\psi).
\end{align}
\end{Lem}
\begin{proof}
Applying \(iD\) to \eqref{psi-thirring-geometric} we find
\begin{align*}
\nabla^2_{\partial_t}\psi-\nabla^2_{\partial_x}\psi=\lambda^2\psi+\epsilon_j\lambda\kappa\langle\psi,e_j\cdot\psi\rangle e_j\cdot\psi+i\kappa D(\epsilon_j\langle\psi,e_j\cdot\psi\rangle e_j\cdot\psi).
\end{align*}
In order to manipulate the last contribution on the right hand side we calculate
\begin{align*}
(\partial_t\cdot\nabla_{\partial_t}-\partial_x\cdot\nabla_{\partial_x})&(\langle\psi,\partial_t\cdot\psi\rangle\partial_t\cdot\psi-\langle\psi,\partial_x\cdot\psi\rangle\partial_x\cdot\psi) \\
=&2\lambda\langle i\partial_x\cdot\partial_t\cdot\psi,\psi\rangle\partial_x\cdot\partial_t\cdot\psi
-\langle\psi,\partial_t\cdot\psi\rangle(\nabla_{\partial_t}\psi+\partial_x\cdot\partial_t\cdot\nabla_{\partial_x}\psi) \\
&+\langle\psi,\partial_x\cdot\psi\rangle(\nabla_{\partial_x}\psi+\partial_x\cdot\partial_t\cdot\nabla_{\partial_t}\psi),
\end{align*}
where we made use of \eqref{identity-tt-xx-psi-thirring-geometric} and \eqref{identity-tx-xt-psi-thirring-geometric}.
Rewriting \eqref{psi-thirring-geometric} as 
\begin{align*}
\nabla_{\partial_x}\psi=\partial_x\cdot\partial_t\cdot\nabla_{\partial_t}\psi+i\lambda\partial_x\cdot\psi+i\kappa\epsilon_j\langle\psi,e_j\cdot\psi\rangle\partial_x\cdot e_j\cdot\psi
\end{align*}
and using the identity
\begin{align*}
\epsilon_j\langle\psi,\partial_t\cdot\psi\rangle&\langle\psi,e_j\cdot\psi\rangle\partial_t\cdot e_j\cdot\psi  
+\epsilon_j\langle\psi,\partial_x\cdot\psi\rangle\langle\psi,e_j\cdot\psi\rangle\partial_x\cdot e_j\cdot\psi \\
=&-|\langle\psi,\partial_t\cdot\psi\rangle|^2\psi-|\langle\psi,\partial_x\cdot\psi\rangle|^2\psi
+2\langle\psi,\partial_t\cdot\psi\rangle\langle\psi,\partial_x\cdot\psi\rangle\partial_x\cdot\partial_t\cdot\psi
\end{align*}
yields the claim.
\end{proof}

\begin{Lem}
Let \(\psi\in\Gamma(\Sigma\R^{1,1})\) be a solution of \eqref{psi-thirring-geometric}. Then the following equality holds
\begin{align}
\label{lemma-thirring-l6psi}
\frac{d}{dt}\int_\R(\frac{1}{3}|\psi|_\beta^6+|\langle\partial_x\cdot\psi,\psi\rangle|^2|\psi|^2_\beta)dx=
4\lambda\int_\R\langle\partial_t\cdot\psi,\psi\rangle\langle\partial_x\cdot\psi,\psi\rangle\langle i\partial_x\cdot\partial_t\cdot\psi,\psi\rangle dx.
\end{align}
\end{Lem}
\begin{proof}
This follows by a direct calculation using \eqref{identity-tt-xx-psi-thirring-geometric} and \eqref{identity-tx-xt-psi-thirring-geometric}.
\end{proof}

\begin{Prop}
Let \(\psi\in\Gamma(\Sigma\R^{1,1})\) be a solution of \eqref{psi-thirring-geometric} with \(\lambda\neq 0\). Then the following inequality holds
\begin{align}
\int_\R(|\nabla_{\partial_t}\psi|^2_\beta+|\nabla_{\partial_x}\psi|^2_\beta)dx\leq Ce^{Ct},
\end{align}
where the constant \(C\) depends on \(\lambda,\kappa\) and the initial data.
\end{Prop}
\begin{proof}
From \eqref{thirring-wave} and a direct calculation we obtain
\begin{align*}
\frac{d}{dt}\frac{1}{2}\int_\R(|\nabla_{\partial_t}\psi|^2_\beta+|\nabla_{\partial_x}\psi|^2_\beta)dx\leq &
\frac{d}{dt}\kappa^2\int_\R(\frac{1}{3}|\psi|_\beta^6+|\langle\partial_x\cdot\psi,\psi\rangle|^2|\psi|^2_\beta)dx
+C\int_\R|\psi|_\beta^3|\nabla\psi|_\beta dx \\
\leq& C\int_\R|\psi|_\beta^6 dx+C\int_\R|\nabla\psi|_\beta^2dx,
\end{align*}
where we used \eqref{lemma-thirring-l6psi} in the last step.
In order to estimate the \(L^6\)-norm of \(\psi\) we make use of the 
Sobolev embedding theorem in one dimension
\begin{align*}
\big(\int_\R|\psi|^6_\beta dx\big)^\frac{1}{2}&\leq C\big(\int_\R\big|\partial_x|\psi|^3_\beta\big|^\frac{2}{3}dx\big)^\frac{3}{2} \\
&\leq C\big(\int_\R|\nabla_{\partial_x}\psi|^\frac{2}{3}_\beta|\psi|^\frac{4}{3}_\beta dx\big)^\frac{3}{2} \\
&\leq C\big(\int_\R|\nabla_{\partial_x}\psi|_\beta^2dx\big)^\frac{1}{2}\int_\R|\psi|_\beta^2dx.
\end{align*}
Since the \(L^2\)-norm of \(\psi\) is conserved for a solution of \eqref{psi-thirring-geometric} we obtain
\begin{align*}
\frac{d}{dt}\frac{1}{2}\int_\R(|\nabla_{\partial_t}\psi|^2_\beta+|\nabla_{\partial_x}\psi|^2_\beta)dx\leq &
C\int_\R(|\nabla_{\partial_t}\psi|^2_\beta+|\nabla_{\partial_x}\psi|^2_\beta)dx
\end{align*}
and the result follows by integration of the differential inequality.
\end{proof}

\begin{Cor}
Let \(\psi\in\Gamma(\Sigma\R^{1,1})\) be a solution of \eqref{psi-thirring-geometric} with \(\lambda\neq 0\).
Then the following estimate holds
\begin{align}
|\psi|_\beta\leq Ce^{Ct},
\end{align}
where the positive constant \(C\) depends on \(\lambda,\kappa\) and the initial data.
\end{Cor}

\begin{Satz}[Existence of a global solution]
\label{existence-thirring}
For any given initial data of the regularity
\begin{align*}
\psi(0,x)=&\psi_0(x)\in H^1(\R,\Sigma\R^{1,1})
\end{align*}
the equation \eqref{psi-thirring-geometric} admits a global weak solution
in \(H^1(\R^{1,1},\Sigma\R^{1,1})\), which is uniquely determined by the initial data.
\end{Satz}

\begin{proof}
The existence of a global solution follows directly since we have a uniform bound on \(\psi\) for the massless case \(\lambda=0\).
Moreover, in the massive case \(\lambda\neq 0\) we have an exponential bound on \(\psi\).
These bounds ensure that the solution cannot blow up and has to exist globally.

To achieve uniqueness let us consider two solutions \(\psi,\xi\) of \eqref{psi-thirring-geometric}
that coincide at \(t=0\). Set \(\eta:=\psi-\xi\). Then \(\eta\) satisfies
\begin{align*}
\partial_t\cdot\nabla_{\partial_t}\eta=\partial_x\cdot\nabla_{\partial_x}\eta-i\lambda\eta
-i\kappa\epsilon_j(\langle\eta,e_j\cdot\psi\rangle e_j\cdot\psi+\langle\xi,e_j\cdot\eta\rangle e_j\cdot\psi
-\langle\xi,e_j\cdot\xi\rangle e_j\cdot\eta).
\end{align*}
Thus, we find
\begin{align*}
\frac{d}{dt}\int_\R|\eta|^2_\beta dx=&2\kappa\int_\R(\epsilon_j\operatorname{Re}\langle\eta,e_j\cdot\psi\rangle\langle ie_j\cdot\psi,\eta\rangle)
+\epsilon_j\operatorname{Re}\langle\xi,e_j\cdot\eta\rangle\langle ie_j\cdot\psi,\eta\rangle) dx \\
\leq& C\int_\R|\eta|^2_\beta(|\psi|^2_\beta+|\chi|_\beta|\psi|_\beta)dx \\
\leq& C\int_\R|\eta|^2_\beta dx,
\end{align*}
where we used the pointwise bound on \(\psi,\xi\) in the last step.
Consequently, we find
\begin{align*}
\int_\R|\eta|^2_\beta dx\leq e^{Ct}\int_\R|\eta|_\beta^2|_{t=0}dx
\end{align*}
such that if \(\psi=\xi\) at \(t=0\) then \(\psi=\xi\) for all times.
\end{proof}

\begin{Bem}
The existence of a global solution for the Thirring model is due to the algebraic structure
of the right hand side of \eqref{psi-thirring-geometric}.
More generally, we could consider
\begin{align*}
iD\psi=\kappa V(\psi)\epsilon_j\langle e_j\cdot\psi,\psi\rangle e_j\cdot\psi,
\end{align*}
where \(V(\psi)\) is supposed to be a real-valued potential. It can again be checked that
\begin{align*}
\Box|\psi|^2_\beta=0 
\end{align*}
such that we get a global bound on \(\psi\).
\end{Bem}

\begin{Bem}
In the physics literature the Thirring model is usually formulated as 
\begin{align*}
S(\psi,\bar\psi)=\int_{\R^{1,1}}(\langle\bar\psi,iD\psi\rangle-\lambda|\psi|_\beta^2-\frac{\kappa}{2}\epsilon_j\langle\bar\psi,e_j\cdot\psi\rangle\langle\bar\psi,e_j\cdot\psi\rangle)d\mu,
\end{align*}
where it is assumed that \(\psi\) and \(\bar\psi\) are independent.
Hence, the critical points consist of two equations
\begin{align*}
iD\psi=\lambda\psi+\kappa\langle\bar\psi,e_j\cdot\psi\rangle e_j\cdot\psi, \qquad iD\bar\psi=\lambda\bar\psi+\kappa\langle\psi,e_j\cdot\bar\psi\rangle e_j\cdot\bar\psi.
\end{align*}
These equations have to be considered as independent as can be checked by a direct calculation.
\end{Bem}

\section{Dirac-wave maps from two-dimensional Minkowski space}
Dirac-wave maps arise as a mathematical version of the supersymmetric nonlinear sigma model
studied in quantum field theory, see for example \cite{MR626710} for the physics background.
The central object of the supersymmetric nonlinear sigma model is an energy functional
that consists of a map between two manifolds and so-called vector spinors.
We want to analyze this model with the methods from geometric analysis, hence
in contrast to the physics literature we will consider standard instead of Grassmann-valued spinors.

In order to define Dirac-wave maps from two-dimensional Minkowski space we choose the following setup. 
Let \((N,h)\) be a compact Riemannian manifold and let \(\phi\colon\R^{1,1}\to N\) be a map.
We consider the pullback of the tangent bundle from the target, which will be denoted by \(\phi^{\ast}TN\).
As discussed in section \ref{subsection-twisted-spinors} we form the twisted spinor bundle \(\Sigma M\otimes\phi^{\ast} TN\).
Sections in \(\Sigma M\otimes\phi^{\ast} TN\) will be called \emph{vector spinors}.

Most of the results that have been obtained in the mathematical literature on the 
supersymmetric nonlinear sigma model consider the case where both domain and target manifolds are Riemannian. 
This study was initiated in \cite{MR2262709}, where the notion of \emph{Dirac-harmonic maps} was introduced.
Dirac-harmonic maps form a semilinear elliptic system for a map between two Riemannian manifolds
and a spinor along that map. For a given Dirac-harmonic map many analytic and geometric results
have been established, as for example the regularity of weak solutions \cite{MR2544729}.
Motivated from the physics literature there exist several extensions of the Dirac-harmonic map system
such as Dirac-harmonic maps with curvature term
\cite{MR3333092}, \cite{MR2370260} and Dirac-harmonic maps with torsion \cite{MR3493217}.

Making use of the Atiyah-Singer index theorem uncoupled solutions to the equations for Dirac-harmonic maps
have been constructed in \cite{MR3070562}. Here, uncoupled refers to the fact that the map part is harmonic.
In addition, several approaches to the existence problem that 
make use of the heat-flow method have been studied in \cite{MR3435758,MR3412386,Jost2017,branding2014evolution,jowi2017} 
and \cite{branding2012evolution}. 

Up to now there is only one reference investigating Dirac-wave maps \cite{MR2138082},
that is critical points of the supersymmetric nonlinear sigma model 
with the domain being two-dimensional Minkowski space.
Expressing the Dirac-wave map system in characteristic coordinates
an existence result for smooth initial data could be obtained.
In this section we will extend the analysis of Dirac-wave maps
and derive an existence result that also includes distributional initial data.
The methods we use here are partly inspired from the analysis of wave maps, see \cite{MR1674843}
for an introduction to the latter.

A problem similar to the one studied in this section, namely the full bosonic string 
from two-dimensional Minkowski space to Riemannian manifolds was treated in \cite{MR3624770}.

The energy functional for Dirac-wave maps is given by
\begin{align}
\label{functional-dwmap-geometric}
S(\phi,\psi)=\frac{1}{2}\int_{\R^{1,1}}(|d\phi|^2+\langle\psi,iD^{\phi^\ast TN}\psi\rangle)d\mu.
\end{align}
Here, \(D^{\phi^\ast TN}\) is the twisted Dirac operator acting on vector spinors. Note that
\(iD^{\phi^\ast TN}\) is self-adjoint with respect to the \(L^2\)-norm such that the energy
functional is real-valued. 

Whenever choosing local coordinates we will use Latin indices to denote coordinates on two-dimensional
Minkowski space and Greek indices to denote coordinates on the target manifold.

For the sake of completeness we will give a short derivation of the critical points of \eqref{functional-dwmap-geometric}.
\begin{Prop}
The Euler-Lagrange equations of \eqref{functional-dwmap-geometric} read
\begin{align}
\label{phi-dwmap}
\tau(\phi)=&\frac{1}{2}R^N(\psi,ie_j\cdot\psi)\epsilon_jd\phi(e_j), \\
\label{psi-dwmap}
D^{\phi^\ast TN}\psi=&0. 
\end{align}
Here, \(R^N\) denotes the curvature tensor of the target \(N\) and \(e_j,j=1,2\) is an pseudo-orthonormal basis of \(T\R^{1,1}\).
\end{Prop}
\begin{proof}
First, we consider a variation of \(\psi\), while keeping \(\phi\) fixed,  satisfying \(\frac{\tilde{\nabla}\psi}{\partial s}\big|_{s=0}=\xi\).
We calculate
\begin{align*}
\frac{d}{ds}\big|_{s=0}S(\phi,\psi)=\frac{1}{2}\int_{\R^{1,1}}(\langle\xi,iD^{\phi^\ast TN}\psi\rangle+\langle\psi,iD^{\phi^\ast TN}\xi\rangle)d\mu
=\int_{\R^{1,1}}\operatorname{Re}\langle\xi,iD^{\phi^\ast TN}\psi\rangle d\mu
\end{align*}
yielding the equation for the vector spinor \(\psi\).
To obtain the equation for the map \(\phi\) we consider a variation of \(\phi\), while keeping \(\psi\) fixed, that is \(\frac{\partial\phi}{\partial s}\big|_{s=0}=\eta\).
It is well-known that
\begin{align*}
\frac{d}{ds}\big|_{s=0}\frac{1}{2}\int_{\R^{1,1}}|d\phi|^2d\mu=-\int_{\R^{1,1}}\langle\tau(\phi),\eta\rangle d\mu,
\end{align*}
where \(\tau(\phi):=\nabla^{\phi^\ast TN}_{e_j}d\phi(e_j)\) is the wave-map operator.
In addition, we find
\begin{align*}
\frac{d}{ds}\big|_{s=0}\frac{1}{2}\int_{\R^{1,1}}\langle\psi,iD^{\phi^\ast TN}\psi\rangle d\mu
=&\frac{1}{2}\int_{\R^{1,1}}\langle\psi,R^N(d\phi(\partial_s),d\phi(e_j))i\epsilon_j e_j\cdot\psi\rangle|_{s=0}d\mu \\
=&\frac{1}{2}\int_{\R^{1,1}}\langle\psi,R^N(\psi,i\epsilon_j e_j\cdot\psi)d\phi(e_j),\eta\rangle d\mu,
\end{align*}
completing the proof.
\end{proof}

Solutions of the system \eqref{phi-dwmap}, \eqref{psi-dwmap} are called \emph{Dirac-wave maps} from \(\R^{1,1}\to N\).

The system \eqref{phi-dwmap}, \eqref{psi-dwmap} can be expanded as 
\begin{align}
\label{el-phi-minkowski}\frac{\nabla}{\partial t}d\phi(\partial_t)-\frac{\nabla}{\partial x}d\phi(\partial_x)&=
\frac{1}{2}R^N(\psi,i\partial_t\cdot\psi)d\phi(\partial_t)
-\frac{1}{2}R^N(\psi,i\partial_x\cdot\psi)d\phi(\partial_x),\\
\label{el-psi-minkowski}\partial_t\cdot\tilde{\nabla}_{\partial_t}\psi&=\partial_x\cdot\tilde{\nabla}_{\partial_x}\psi.
\end{align}

Choosing local coordinates on the target \(N\), the Euler-Lagrange equations acquire the form \((\alpha=1,\ldots,\dim N)\)
\begin{align*}
\frac{\partial^2\phi^\alpha}{\partial t^2}-\frac{\partial^2\phi^\alpha}{\partial x^2}
+\Gamma^\alpha_{\beta\gamma}\big(\frac{\partial\phi^\beta}{\partial t}\frac{\partial\phi^\gamma}{\partial t}
-\frac{\partial\phi^\beta}{\partial x}\frac{\partial\phi^\gamma}{\partial x}\big)=&
R^\alpha_{~\beta\gamma\delta}(\langle\psi^\gamma,i\partial_t\cdot\psi^\delta\rangle\frac{\partial\phi^\beta}{\partial t}
-\langle\psi^\gamma,i\partial_x\cdot\psi^\delta\rangle\frac{\partial\phi^\beta}{\partial x}),\\
D\psi^\alpha=&-\Gamma^\alpha_{\beta\gamma}\frac{\partial\phi^\beta}{\partial x_j}\epsilon_je_j\cdot\psi^\gamma,
\end{align*}
where \(\Gamma^\alpha_{\beta\gamma}\) are the Christoffel symbols and \(R^\alpha_{~\beta\gamma\delta}\) the components of the
curvature tensor on the target \(N\).

In order to treat a weak version of the system \eqref{phi-dwmap}, \eqref{psi-dwmap}
it will be necessary to embed \(N\) isometrically into some \(\R^q\) making use of the Nash embedding theorem. 
Then we have \(\phi\colon\R^{1,1}\to\R^q\) and \(\psi\in\Gamma(\Sigma\R^{1,1}\otimes\R^q)\).
In this case the equations for Dirac-wave maps acquire the form
\begin{align}
\label{phi-dwmap-extrinsic}
\frac{\partial^2\phi}{\partial t^2}-\frac{\partial^2\phi}{\partial x^2}=&\sff(d\phi,d\phi)
+P(\sff(\psi,d\phi(\partial_t)),i\partial_t\cdot\psi)-P(\sff(\psi,d\phi(\partial_x)),i\partial_x\cdot\psi), \\
\label{psi-dwmap-extrinsic}
\partial_t\cdot\nabla_{\partial_t}\psi-\partial_x\cdot\nabla_{\partial_x}\psi
=&\sff(\frac{\partial\phi}{\partial t},\partial_t\cdot\psi)-\sff(\frac{\partial\phi}{\partial x},\partial_x\cdot\psi),
\end{align}
where \(\sff\) denotes the second fundamental form of the embedding and
\(P\) the shape operator defined by
\begin{align*}
\langle P(\xi,X),Y)\rangle=\langle\sff(X,Y),\xi\rangle
\end{align*}
for \(X,Y\in\Gamma(TN)\) and \(\xi\in T^\perp N\).
For solutions of \eqref{psi-dwmap} we obtain the following identities
\begin{align}
\label{identity-tt-xx-psi-dwmap} \frac{\partial}{\partial t}\langle\partial_t\cdot\psi,\psi\rangle-\frac{\partial}{\partial x}\langle\partial_x\cdot\psi,\psi\rangle
=&\langle D^{\phi^\ast TN}\psi,\psi\rangle+\langle\psi, D^{\phi^\ast TN}\psi\rangle=0, \\
\label{identity-tx-xt-psi-dwmap}\frac{\partial}{\partial t}\langle\partial_x\cdot\psi,\psi\rangle-\frac{\partial}{\partial x}\langle\partial_t\cdot\psi,\psi\rangle
=&\langle(\partial_t\cdot\tilde{\nabla}_{\partial_x}-\partial_x\cdot\tilde{\nabla}_{\partial_t})\psi,\psi\rangle
-\langle\psi,(\partial_t\cdot\tilde{\nabla}_{\partial_x}+\partial_x\cdot\tilde{\nabla}_{\partial_t})\psi\rangle=0.
\end{align}
Moreover, for a solution of \eqref{psi-dwmap} the following wave type equation holds
\begin{align*}
\tilde{\nabla}_{\partial_t}^2\psi-\tilde{\nabla}_{\partial_x}^2\psi=-\partial_t\cdot\partial_x\cdot R^N(d\phi(\partial_t),d\phi(\partial_x))\psi
\end{align*}
due to the Weitzenboeck formula \eqref{weitzenboeck-twisted} with \(F=\phi^\ast TN\).

\begin{Bem}
In the physics literature \cite{MR626710} the energy functional for the supersymmetric nonlinear sigma model is defined as follows
\begin{align*}
S(\phi,\psi,\bar\psi)=\frac{1}{2}\int_{\R^{1,1}}(|d\phi|^2+\langle\bar{\psi},iD^{\phi^\ast TN}\psi\rangle)d\mu,
\end{align*}
where \(\bar\psi:=\partial_t\cdot\psi\). Treating \(\psi,\bar\psi\) as independent fields, we obtain the Euler-Lagrange equations
\begin{align*}
\tau(\phi)=&\frac{1}{2}R^N(\bar\psi,ie_j\cdot\psi)\epsilon_jd\phi(e_j), \qquad D^{\phi^\ast TN}\psi=0,\qquad D^{\phi^\ast TN}\bar\psi=0.
\end{align*}
Note that the spinors \(\psi\) and \(\bar\psi\) have to be considered as independent since the two equations
\(D^{\phi^\ast TN}\psi=0, ~~ D^{\phi^\ast TN}(\partial_t\cdot\psi)=0\) are not compatible.
In addition, the equation for the map \(\phi\) would acquire the form 
\begin{align*}
\tau(\phi)=&-\underbrace{\frac{1}{2}R^N(\psi,i\partial_t\cdot\partial_x\cdot\psi)d\phi(\partial_x)}_{=0},
\end{align*}
when inserting \(\bar\psi:=\partial_t\cdot\psi\).
It turns out that the curvature term on the right hand side vanishes due to symmetry reasons.
In the physics literature one usually considers anticommuting spinors such that this term does not vanish.
\end{Bem}
In the following we will only analyze the system \eqref{phi-dwmap}, \eqref{psi-dwmap}.

Let us demonstrate how to construct an explicit solution to the Euler-Lagrange equations \eqref{phi-dwmap}, \eqref{psi-dwmap},
where we follow the ideas used for the construction of an explicit solution to the Dirac-harmonic map system in \cite[Proposition 2.2]{MR2262709}.
To this end we recall the following 
\begin{Dfn}
Let \(M\) be a \(n\)-dimensional Lorentzian spin manifold.
A spinor \(\psi\in\Gamma(\Sigma M)\) is called \emph{twistor spinor}
if it satisfies 
\begin{align}
P_X\psi:=\nabla^{\Sigma M}_X\psi+\frac{1}{n}X\cdot D\psi=0
\end{align}
for all vector fields \(X\).
\end{Dfn}

In two-dimensional Minkowski space twistor spinors are of the form
\begin{align*}
\psi(x)=\psi_1+x\cdot\psi_2,
\end{align*}
where \(\psi_1,\psi_2\) are constant spinors \cite{MR1703175}.

\begin{Prop}
Let \(\phi\colon\R^{1,1}\to N\) be a wave map, that is a solution of \(\tau(\phi)=0\).
We set
\begin{align*}
\psi:=\epsilon_je_j\cdot\chi\otimes d\phi(e_j),
\end{align*}
where \(\chi\) is a twistor spinor. 
Then the pair \((\phi,\psi)\) is a Dirac-wave map, that is uncoupled:
\begin{align*}
\tau(\phi)=0=\frac{1}{2}\epsilon_jR^N(\psi,ie_j\cdot\psi)d\phi(e_j),\qquad D^{\phi^\ast TN}\psi=0.
\end{align*}
\end{Prop}
\begin{proof}
First, we check that the equation for the vector spinor \(\psi\) is satisfied.
To this end, we calculate
\begin{align*}
D^{\phi^\ast TN}\psi=&\big(\partial_t\cdot\tilde{\nabla}_{\partial_t}-\partial_x\cdot\tilde{\nabla}_{\partial_x}\big)
\big(\partial_t\cdot\chi\otimes d\phi(\partial_t)-\partial_x\cdot\chi\otimes d\phi(\partial_x)\big) \\
=&-\big(\nabla^{\Sigma M}_{\partial_t}\chi+\partial_x\cdot\partial_t\cdot\nabla^{\Sigma M}_{\partial_x}\chi\big)\otimes d\phi(\partial_t)
+\big(\nabla^{\Sigma M}_{\partial_x}\chi-\partial_t\cdot\partial_x\cdot\nabla^{\Sigma M}_{\partial_t}\chi\big)\otimes d\phi(\partial_x)\\
&-\chi\otimes\tau(\phi)-\partial_t\cdot\partial_x\cdot\chi\otimes(\frac{\nabla}{\partial t}d\phi(\partial_x)-\frac{\nabla}{\partial x}d\phi(\partial_t))\\
=&0,
\end{align*}
the first two terms vanish since \(\chi\) is a twistor spinor by assumption.
As a second step we check that the curvature term on the right hand side of \eqref{phi-dwmap} vanishes.
Using the local expression of \eqref{phi-dwmap}
\begin{align*}
\epsilon_jR^N(\psi,ie_j\cdot\psi)d\phi(e_j)=&R^\alpha_{~\beta\gamma\delta}\frac{\partial}{\partial y^\alpha}
(\langle\psi^\gamma,i\partial_t\cdot\psi^\delta\rangle\frac{\partial\phi^\beta}{\partial t}
-\langle\psi^\gamma,i\partial_x\cdot\psi^\delta\rangle\frac{\partial\phi^\beta}{\partial x})
\end{align*}
we find (the second term vanishes for the same reason)
\begin{align*}
R^\alpha_{~\beta\gamma\delta}\langle\psi^\gamma,i\partial_t\cdot\psi^\delta\rangle=&
-R^\alpha_{~\beta\gamma\delta}\langle\chi,i\partial_t\cdot\chi\rangle
\big(\frac{\partial\phi^\gamma}{\partial t}\frac{\partial\phi^\delta}{\partial t}
+\frac{\partial\phi^\gamma}{\partial x}\frac{\partial\phi^\delta}{\partial x}\big) \\
&+R^\alpha_{~\beta\gamma\delta}\langle\chi,i\partial_x\cdot\chi\rangle
\big(\frac{\partial\phi^\gamma}{\partial t}\frac{\partial\phi^\delta}{\partial x}
+\frac{\partial\phi^\gamma}{\partial x}\frac{\partial\phi^\delta}{\partial t}\big)=0
\end{align*}
due to the symmetries of the curvature tensor on \(N\).
\end{proof}

\subsection{Conserved Energies}
In this subsection we give several conserved energies for solutions of the system \eqref{phi-dwmap}, \eqref{psi-dwmap}.
By integrating \eqref{identity-tt-xx-psi-dwmap} we directly get that
\begin{align*}
E_1(t)=\int_\R|\psi|^2_\beta dx
\end{align*}
is conserved for a solution of \eqref{psi-dwmap}.
Again, it is straightforward to also control higher \(L^p\) norms of \(\psi\) as discussed in \ref{subsection-twisted-spinors}.

Moreover, for a solution of \eqref{psi-dwmap} we have
\begin{align*}
|\tilde{\nabla}_{\partial_t}\psi|^2_\beta=|\tilde{\nabla}_{\partial_x}\psi|^2_\beta.
\end{align*}

\begin{Prop}
Let \((\phi,\psi)\colon\R^{1,1}\to N\) be a Dirac-wave map.
Then the energy 
\begin{align}
\label{dwmap-energy-intrinsic}
E(t):=\frac{1}{2}\int_\R\big(|d\phi(\partial_t)|^2+|d\phi(\partial_x)|^2-\langle\psi,i\partial_t\cdot\tilde{\nabla}_{\partial_t}\psi\rangle\big)dx
\end{align}
is conserved.
\end{Prop}
\begin{proof}
We calculate
\begin{align*}
\frac{d}{dt}\frac{1}{2}\int_\R(|d\phi(\partial_t)|^2+|d\phi(\partial_x)|^2)dx
=&\int_\R\langle d\phi(\partial_t),\tau(\phi)\rangle dx \\
=&-\frac{1}{2}\int_\R\langle d\phi(\partial_t),R^N(\psi,i\partial_x\cdot\psi)d\phi(\partial_x)\rangle dx.
\end{align*}
Differentiating \eqref{psi-dwmap} with respect to \(t\) we obtain the following identity
\begin{align*}
\partial_t\cdot\tilde\nabla^2_{\partial_t}\psi=\partial_x\cdot\tilde\nabla_{\partial_t}\tilde\nabla_{\partial_x}\psi
=\partial_x\cdot R^N(d\phi(\partial_t),d\phi(\partial_x))\psi+\partial_x\cdot\tilde{\nabla}_{\partial_x}\tilde{\nabla}_{\partial_t}\psi.
\end{align*}
Thus, we find
\begin{align*}
\frac{d}{dt}\frac{1}{2}\int_\R\langle\psi,i\partial_t\cdot\tilde{\nabla}_{\partial_t}\psi\rangle dx
=\frac{1}{2}\int_\R&(\langle\underbrace{(\partial_t\cdot\tilde{\nabla}_{\partial_t}-\partial_x\cdot\tilde{\nabla}_{\partial_x})\psi}_{=0},i\tilde{\nabla}_{\partial_t}\psi\rangle \\
&+\langle\psi,i\partial_x\cdot R^N(d\phi(\partial_t),d\phi(\partial_x))\psi\rangle)dx,
\end{align*}
where we used integration by parts in the last step. The assertion then follows by adding up both contributions.
\end{proof}

Having gained control over \(\psi\) we now establish a bound on the derivatives of \(\phi\). To this end we set 
\begin{align*}
e(\phi):=\frac{1}{2}(|d\phi(\partial_t)|^2+|d\phi(\partial_x)|^2).
\end{align*}

Then we obtain the following
\begin{Prop}
Let \((\phi,\psi)\colon\R^{1,1}\to N\) be a Dirac-wave map. Then the following formula holds
\begin{align}
\label{dwmap-box-energy-phi}
\Box e(\phi)=\frac{\partial^2}{\partial t^2}\langle\tilde{\nabla}_{\partial_x}\psi,i\partial_x\cdot\psi\rangle
-\frac{\partial^2}{\partial x^2}\langle\tilde{\nabla}_{\partial_t}\psi,i\partial_t\cdot\psi\rangle,
\end{align}
where \(\Box:=\frac{\partial^2}{\partial t^2}-\frac{\partial^2}{\partial x^2}\).
\end{Prop}

\begin{proof}
Testing \eqref{el-phi-minkowski} with \(d\phi(\partial_t)\) and \(d\phi(\partial_x)\) we obtain the two equations
\begin{align*}
\frac{\partial}{\partial t}e(\phi)-\frac{\partial}{\partial x}\langle d\phi(\partial_t),d\phi(\partial_x)\rangle
&=-\langle R^N(\psi,i\partial_x\cdot\psi)d\phi(\partial_x),d\phi(\partial_t)\rangle, \\
-\frac{\partial}{\partial x}e(\phi)+\frac{\partial}{\partial t}\langle d\phi(\partial_t),d\phi(\partial_x)\rangle
&=\langle R^N(\psi,i\partial_t\cdot\psi)d\phi(\partial_t),d\phi(\partial_x)\rangle.
\end{align*}
Differentiating the first equation with respect to \(t\), the second one with respect to \(x\) and adding up both contributions we find
\begin{align*}
\Box e(\phi)=\frac{\partial}{\partial x}\langle R^N(\psi,i\partial_t\cdot\psi)d\phi(\partial_t),d\phi(\partial_x)\rangle
-\frac{\partial}{\partial t}\langle R^N(\psi,i\partial_x\cdot\psi)d\phi(\partial_x),d\phi(\partial_t)\rangle.
\end{align*}
We may rewrite the right-hand side as follows
\begin{align*}
&\frac{\partial}{\partial x}\langle R^N(\psi,i\partial_t\cdot\psi)d\phi(\partial_t),d\phi(\partial_x)\rangle
-\frac{\partial}{\partial t}\langle R^N(\psi,i\partial_x\cdot\psi)d\phi(\partial_x),d\phi(\partial_t)\rangle \\
=&\frac{\partial}{\partial x}\langle\tilde{\nabla}_{\partial_t}\tilde{\nabla}_{\partial_x}\psi,i\partial_t\cdot\psi\rangle
-\frac{\partial}{\partial x}\langle\tilde{\nabla}_{\partial_x}\tilde{\nabla}_{\partial_t}\psi,i\partial_t\cdot\psi\rangle 
-\frac{\partial}{\partial t}\langle\tilde{\nabla}_{\partial_x}\tilde{\nabla}_{\partial_t}\psi,i\partial_x\cdot\psi\rangle
+\frac{\partial}{\partial t}\langle\tilde{\nabla}_{\partial_t}\tilde{\nabla}_{\partial_x}\psi,i\partial_x\cdot\psi\rangle \\
=&\frac{\partial^2}{\partial x\partial t}\langle\tilde{\nabla}_{\partial_x}\psi,i\partial_t\cdot\psi\rangle 
-\frac{\partial}{\partial x}\langle\tilde{\nabla}_{\partial_x}\psi,i\partial_t\cdot\tilde{\nabla}_{\partial_t}\psi\rangle
-\frac{\partial^2}{\partial x^2}\langle\tilde{\nabla}_{\partial_t}\psi,i\partial_t\cdot\psi\rangle 
+\frac{\partial}{\partial x}\langle\tilde{\nabla}_{\partial_t}\psi,i\partial_t\cdot\tilde{\nabla}_{\partial_x}\psi\rangle \\
&-\frac{\partial^2}{\partial x\partial t}\langle\tilde{\nabla}_{\partial_t}\psi,i\partial_x\cdot\psi\rangle 
+\frac{\partial}{\partial t}\langle\tilde{\nabla}_{\partial_t}\psi,i\partial_x\cdot\tilde{\nabla}_{\partial_x}\psi\rangle
+\frac{\partial^2}{\partial t^2}\langle\tilde{\nabla}_{\partial_x}\psi,i\partial_x\cdot\psi\rangle 
-\frac{\partial}{\partial t}\langle\tilde{\nabla}_{\partial_x}\psi,i\partial_x\cdot\tilde{\nabla}_{\partial_t}\psi\rangle \\
=&\frac{\partial^2}{\partial x\partial t}\langle\underbrace{(\partial_t\cdot\tilde{\nabla}_{\partial_x}-\partial_x\cdot\tilde{\nabla}_{\partial_t})\psi}_{=0},i\psi\rangle 
+\frac{\partial}{\partial x}\big(\underbrace{\langle\tilde{\nabla}_{\partial_t}\psi,i\partial_t\cdot\tilde{\nabla}_{\partial_x}\psi\rangle 
-\langle\tilde{\nabla}_{\partial_x}\psi,i\partial_t\cdot\tilde{\nabla}_{\partial_t}\psi\rangle}_{=0}\big) \\
&+\frac{\partial}{\partial t}\big(\underbrace{\langle\tilde{\nabla}_{\partial_t}\psi,i\partial_x\cdot\tilde{\nabla}_{\partial_x}\psi\rangle
-\langle\tilde{\nabla}_{\partial_x}\psi,i\partial_x\cdot\tilde{\nabla}_{\partial_t}\psi\rangle}_{=0}\big) 
+\frac{\partial^2}{\partial t^2}\langle\tilde{\nabla}_{\partial_x}\psi,i\partial_x\cdot\psi\rangle
-\frac{\partial^2}{\partial x^2}\langle\tilde{\nabla}_{\partial_t}\psi,i\partial_t\cdot\psi\rangle \\
=&\frac{\partial^2}{\partial t^2}\langle\tilde{\nabla}_{\partial_x}\psi,i\partial_x\cdot\psi\rangle
-\frac{\partial^2}{\partial x^2}\langle\tilde{\nabla}_{\partial_t}\psi,i\partial_t\cdot\psi\rangle,
\end{align*}
where we used that \(\psi\) is a solution of \eqref{el-psi-minkowski} several times completing the proof.
\end{proof}

\begin{Bem}
The conserved energies that we have presented so far all reflect the hyperbolic nature of the 
Dirac-wave map system \eqref{phi-dwmap}, \eqref{psi-dwmap}. In addition, as in the Riemannian case,
we also have the energy-momentum tensor, which is conserved for a Dirac-wave map. More precisely,
the symmetric 2-tensor \(T_{ij}\) defined by
\begin{align*}
T_{ij}:=2\langle d\phi(e_i),d\phi(e_j)\rangle-g_{ij}|d\phi|^2-\frac{1}{2}\langle\psi,(e_i\cdot\tilde{\nabla}_{e_j}+e_j\cdot\tilde{\nabla}_{e_i})\psi\rangle
\end{align*}
is divergence free for a solution of \eqref{phi-dwmap}, \eqref{psi-dwmap}.
\end{Bem}

\subsection{An existence result via energy methods}
In this subsection we will derive an existence result for the Cauchy problem associated to \eqref{phi-dwmap}, \eqref{psi-dwmap}.
Since we want to be able to treat initial data of low regularity, we have to use the extrinsic version of the 
Dirac-wave map system \eqref{phi-dwmap-extrinsic}, \eqref{psi-dwmap-extrinsic}.

Making use of the conserved energies from the last subsection we find
\begin{Prop}
Let \((\phi,\psi)\colon\R^{1,1}\to\R^q\) be a Dirac-wave map. Then the following
energy is conserved
\begin{align}
\label{dwmap-conserved-energy-1}
E_{DW}(t):=\frac{1}{2}\int_\R\big(|\frac{\partial\phi}{\partial t}\big|^2+|\frac{\partial\phi}{\partial x}\big|^2
-\langle\psi,i\partial_t\cdot\nabla_{\partial_t}\psi\rangle\big)dx.
\end{align}
\end{Prop}
\begin{proof}
By \eqref{psi-dwmap-extrinsic} we have a conserved energy for the intrinsic version of the Dirac-wave map system.
To obtain the conserved energy for the extrinsic version we consider the isometric embedding \(\iota\colon N\to\R^q\).
Since \(\iota\) is an isometry we may apply its differential to all terms in \eqref{dwmap-conserved-energy-1},
which gives the statement.
\end{proof}

\begin{Cor}
Let \((\phi,\psi)\colon\R^{1,1}\to\R^q\) be a Dirac-wave map. Then the energy
\begin{align}
\label{dwmap-conserved-energy-2}
E(t)=\int_\R\big(\big|\frac{\partial}{\partial x}(e(\phi)-\langle i\partial_t\cdot\nabla_{\partial_t}\psi,\psi\rangle)\big|^2
+|\frac{\partial}{\partial t}(e(\phi)-\langle i\partial_t\cdot\nabla_{\partial_t}\psi,\psi\rangle)\big|^2\big)dx
\end{align}
is conserved.
\end{Cor}
\begin{proof}
We obtain a conserved energy from \eqref{dwmap-box-energy-phi}. Applying the isometric embedding \(\iota\) again yields the claim.
\end{proof}

\begin{Prop}
Let \(\psi\) be a solution of \eqref{psi-dwmap}. Then the following pointwise bound holds
\begin{align}
\label{psi-bound-dwmap}
|\psi|_\beta\leq C,
\end{align}
where the positive constant \(C\) depends on \(\psi_0\).
\end{Prop}
\begin{proof}
This follows from Proposition \ref{prop-psi-kernel-twisted-a} 
with \(F=\phi^{\ast}TN\) and the Sobolev embedding \(H^1\hookrightarrow L^\infty\).
\end{proof}

\begin{Lem}
Let \((\phi,\psi)\colon\R^{1,1}\to\R^q\) be a Dirac-wave map.
Then the following estimate holds
\begin{align}
\label{dwmap-inequality-phi-psi}
|\frac{\partial\phi}{\partial t}\big|^2+|\frac{\partial\phi}{\partial x}\big|^2\leq C(1+|\nabla\psi|_\beta),
\end{align}
where the positive constant \(C\) depends on \(E_{DW}(0)\) and \(\psi_0\).
\end{Lem}
\begin{proof}
From the conserved energy \eqref{dwmap-conserved-energy-2} and the Sobolev embedding \(H^1\hookrightarrow L^\infty\)
we obtain the following bound
\begin{align*}
|e(\phi)-\langle i\partial_t\cdot\nabla_{\partial_t}\psi,\psi\rangle)|\leq E_{DW}(0).
\end{align*}
Using the pointwise bound on \(\psi\) we obtain the claim.
\end{proof}

In order to obtain control over the derivatives of \(\psi\) we turn \eqref{psi-dwmap-extrinsic} into a wave-type equation.
Applying the Dirac operator \(D\) on both sides of \eqref{psi-dwmap-extrinsic} we obtain the wave-type equation
\begin{align}
\nabla^2_{\partial_t}\psi-\nabla^2_{\partial_x}\psi 
=&\sff(\Box\phi,\psi)+\sff(\frac{\partial\phi}{\partial t},\nabla_{\partial_t}\psi+\partial_x\cdot\partial_t\cdot\nabla_{\partial_x}\psi) 
\nonumber+\sff(\frac{\partial\phi}{\partial x},-\nabla_{\partial_x}\psi+\partial_t\cdot\partial_x\cdot\nabla_{\partial_t}\psi) \\
\nonumber&+(\nabla_{d\phi(\partial_t)}\sff)(\frac{\partial\phi}{\partial t},\psi)
+(\nabla_{d\phi(\partial_t)}\sff)(\frac{\partial\phi}{\partial x},\partial_t\cdot\partial_x\cdot\psi) \\
&-(\nabla_{d\phi(\partial_x)}\sff)(\frac{\partial\phi}{\partial x},\partial_t\cdot\partial_x\cdot\psi)
+(\nabla_{d\phi(\partial_x)}\sff)(\frac{\partial\phi}{\partial t},\psi).
\label{dwmap-wave-equation-spinor-rq}
\end{align}

We set
\begin{align*}
E_{(\psi,1,2)}(t):=\frac{1}{2}\int_\R(|\nabla_{\partial_t}\psi|^2_\beta+|\nabla_{\partial_x}\psi|^2_\beta)dx.
\end{align*}

The energy \(E_{(\psi,1,2)}(t)\) satisfies the following differential inequality
\begin{Lem}
\label{lemma-dwmap-psi12}
Let \((\phi,\psi)\colon\R^{1,1}\to\R^q\) be a Dirac-wave map.
Then the following inequality holds
\begin{align}
\frac{d}{dt}E_{(\psi,1,2)}(t)\leq C\big(E_{(\psi,1,2)}(t)+(E_{(\psi,1,2)}(t))^\frac{1}{2}\big),
\end{align}
where the positive constant \(C\) depends on \(N\) and the initial data.
\end{Lem}
\begin{proof}
Making use of \eqref{dwmap-wave-equation-spinor-rq} a direct calculation yields
\begin{align*}
\frac{d}{dt}E_{(\psi,1,2)}(t)=\int_\R&\big(\langle\partial_t\cdot\nabla_{\partial_t}\psi,\sff(\Box\phi,\psi)\rangle 
+\langle\partial_t\cdot\nabla_{\partial_t}\psi,\sff(\frac{\partial\phi}{\partial t},\nabla_{\partial_t}\psi+\partial_x\cdot\partial_t\cdot\nabla_{\partial_x}\psi)\rangle \\
&+\langle\partial_t\cdot\nabla_{\partial_t}\psi,\sff(\frac{\partial\phi}{\partial x},-\nabla_{\partial_x}\psi+\partial_t\cdot\partial_x\cdot\nabla_{\partial_t}\psi)\rangle 
+\langle\partial_t\cdot\nabla_{\partial_t}\psi,(\nabla_{d\phi(\partial_t)}\sff)(\frac{\partial\phi}{\partial t},\psi)\rangle \\
&+\langle\partial_t\cdot\nabla_{\partial_t}\psi,(\nabla_{d\phi(\partial_t)}\sff)(\frac{\partial\phi}{\partial x},\partial_t\cdot\partial_x\cdot\psi)\rangle 
+\langle\partial_t\cdot\nabla_{\partial_t}\psi,(\nabla_{d\phi(\partial_x)}\sff)(\frac{\partial\phi}{\partial t},\psi)\rangle \\
&-\langle\partial_t\cdot\nabla_{\partial_t}\psi,(\nabla_{d\phi(\partial_x)}\sff)(\frac{\partial\phi}{\partial x},\partial_t\cdot\partial_x\cdot\psi)\rangle  \big)dx.
\end{align*}
Note that the second and the third term on the right hand side can be rewritten as 
\begin{align*}
\langle\partial_t\cdot\nabla_{\partial_t}\psi,\sff(\frac{\partial\phi}{\partial t},\nabla_{\partial_t}\psi+\partial_x\cdot\partial_t\cdot\nabla_{\partial_x}\psi)\rangle
=-\langle\partial_t\cdot\psi,(\nabla_{d\phi(\partial_t)}\sff)(\frac{\partial\phi}{\partial t},\nabla_{\partial_t}\psi+\partial_x\cdot\partial_t\cdot\nabla_{\partial_x}\psi)\rangle
\end{align*}
since \(\psi\perp\sff\). Consequently, we get the following inequality
\begin{align*}
\frac{d}{dt}E_{(\psi,1,2)}(t)
\leq &C\int_\R(|\psi|_\beta|\nabla\psi|_\beta|\Box\phi|+|\psi|_\beta|\nabla\psi|_\beta|d\phi|^2)dx \\
\leq &C\int_\R(|\psi|_\beta|\nabla\psi|_\beta|d\phi|^2 +|\psi|^5_\beta|\nabla\psi|_\beta)dx \\
\leq &C\big(\int_\R|\nabla\psi|^2_\beta dx+\big(\int_\R|\nabla\psi|^2_\beta dx\big)^\frac{1}{2}\big),
\end{align*}
where we used \eqref{dwmap-inequality-phi-psi} in the last step, completing the proof.
\end{proof}

\begin{Cor}
Let \((\phi,\psi)\colon\R^{1,1}\to\R^q\) be a Dirac-wave map.
Then the following estimate holds
\begin{align*}
E_{(\psi,1,2)}(t)\leq Ce^{Ct},
\end{align*}
where the positive constant \(C\) depends on \(N\) and the initial data.
\end{Cor}
\begin{proof}
This follows from the last lemma and the Gronwall inequality.
\end{proof}

\begin{Cor}
Let \((\phi,\psi)\colon\R^{1,1}\to\R^q\) be a Dirac-wave map.
Then the following estimate holds
\begin{align*}
\int_\R\big(|\frac{\partial\phi}{\partial t}\big|^2+|\frac{\partial\phi}{\partial x}\big|^2\big)dx\leq C(1+e^{Ct}),
\end{align*}
where the positive constant \(C\) depends on \(N\) and the initial data.
\end{Cor}
\begin{proof}
From the conserved energy \eqref{dwmap-conserved-energy-1} we obtain
\begin{align*}
\int_\R\big(|\frac{\partial\phi}{\partial t}\big|^2+|\frac{\partial\phi}{\partial x}\big|^2\big)dx
\leq E_{DW}(0)+C\big(\int_\R|\psi|^2_\beta dx\big)^\frac{1}{2}\big(\int_\R|\nabla\psi|^2_\beta)dx\big)^\frac{1}{2} \leq E_{DW}(0)+Ce^{Ct},
\end{align*}
yielding the result.
\end{proof}

It turns out that we need to gain control over the \(L^4\)-norm of \(\nabla\psi\). 
To this end let us recall the following fact.
\begin{Bem}
Suppose that \(f\colon\R^{1,1}\to\R\) is a solution of the scalar wave equation, that is \(\Box f=0\).
Then we have the following conservation law
\begin{align*}
\frac{d}{dt}\int_\R\big(\big|\frac{\partial f}{\partial t}\big|^4+\big|\frac{\partial f}{\partial x}\big|^4
+6\big|\frac{\partial f}{\partial t}\big|^2\big|\frac{\partial f}{\partial x}\big|^2\big)dx=0.
\end{align*}
\end{Bem}

Motivated from the conserved energy for solutions of the scalar wave equation we define
\begin{align*}
E_{(\psi,1,4)}(t):=\int_\R(|\nabla_{\partial_t}\psi|^4_\beta+|\nabla_{\partial_x}\psi|^4_\beta
+2|\nabla_{\partial_x}\psi|^2_\beta|\nabla_{\partial_t}\psi|^2_\beta+4|\langle\partial_t\cdot\nabla_{\partial_x}\psi,\nabla_{\partial_t}\psi\rangle|^2)dx.
\end{align*}
It can be checked by a direct calculation that \(E_{(\psi,1,4)}(t)\) is conserved when \(\psi\in\Gamma(\Sigma\R^{1,1})\)
is a solution of \(\nabla^2_{\partial_t}\psi=\nabla^2_{\partial_x}\psi\).

\begin{Lem}
Let \((\phi,\psi)\colon\R^{1,1}\to\R^q\) be a Dirac-wave map.
Then the following inequality holds
\begin{align}
\frac{d}{dt}E_{(\psi,1,4)}(t)\leq C\big(E_{(\psi,1,4)}(t)+(E_{(\psi,1,4)}(t))^\frac{1}{2}\big),
\end{align}
where the positive constant \(C\) depends on \(N\) and the initial data.
\end{Lem}
\begin{proof}
Let \(\psi\) be a solution of \(\nabla^2_{\partial t}\psi-\nabla^2_{\partial x}\psi=f\).
Then a direct, but lengthy calculation yields
\begin{align*}
\frac{d}{dt}E_{(\psi,1,4)}(t)=
\int_\R(&4|\nabla_{\partial_t}\psi|^2_\beta\langle\partial_t\cdot\nabla_{\partial_t}\psi,f\rangle-4|\nabla_{\partial_x}\psi|^2_\beta\langle\partial_t\cdot\nabla_{\partial_t}\psi,f\rangle \\
&+8\langle\partial_t\cdot\nabla_{\partial_x}\psi,f\rangle\langle\partial_t\cdot\nabla_{\partial_t}\psi,\nabla_{\partial_x}\psi\rangle)dx.
\end{align*}
At this point we choose
\begin{align*}
f=&\sff(\Box\phi,\psi)+\sff(\frac{\partial\phi}{\partial t},\nabla_{\partial_t}\psi+\partial_x\cdot\partial_t\cdot\nabla_{\partial_x}\psi) 
+\sff(\frac{\partial\phi}{\partial x},-\nabla_{\partial_x}\psi+\partial_t\cdot\partial_x\cdot\nabla_{\partial_t}\psi) \\
&+(\nabla_{d\phi(\partial_t)}\sff)(\frac{\partial\phi}{\partial t},\psi)
+(\nabla_{d\phi(\partial_t)}\sff)(\frac{\partial\phi}{\partial x},\partial_t\cdot\partial_x\cdot\psi) \\
&-(\nabla_{d\phi(\partial_x)}\sff)(\frac{\partial\phi}{\partial x},\partial_t\cdot\partial_x\cdot\psi)
+(\nabla_{d\phi(\partial_x)}\sff)(\frac{\partial\phi}{\partial t},\psi).
\end{align*}
By the same reasoning as in the proof of Lemma \ref{lemma-dwmap-psi12} we find
\begin{align*}
\frac{d}{dt}E_{(\psi,1,4)}(t)\leq &C\int_\R(|\psi|_\beta||\nabla\psi|^3_\beta|\Box\phi|+|\psi|_\beta|\nabla\psi|^3_\beta|d\phi|^2)dx \\
\leq &C\int_\R|\nabla\psi|^4_\beta dx+C\big(\int_\R|\nabla\psi|_\beta^4dx\big)^\frac{1}{2},
\end{align*}
which completes the proof.
\end{proof}

\begin{Cor}
Let \((\phi,\psi)\colon\R^{1,1}\to\R^q\) be a Dirac-wave map.
Then the following estimate holds
\begin{align}
E_{(\psi,1,4)}(t)\leq Ce^{Ct},
\end{align}
where the positive constant \(C\) depends on \(N\) and the initial data.
\end{Cor}
\begin{proof}
This follows from the last lemma and the Gronwall inequality.
\end{proof}

After having gained control over the derivatives of \(\psi\) we 
now control the second derivative of \(\phi\). To this end we set
\begin{align*}
E_{(\phi,2,2)}(t):=\frac{1}{2}\int_\R\big(\big|\frac{\partial^2\phi}{\partial x^2}\big|^2+\big|\frac{\partial^2\phi}{\partial x\partial t}\big|^2\big)dx.
\end{align*}

\begin{Prop}
Let \((\phi,\psi)\colon\R^{1,1}\to\R^q\) be a Dirac-wave map.
Then the following inequality holds
\begin{align}
\frac{d}{dt}E_{(\phi,2,2)}(t)\leq C\big(E_{(\phi,2,2)}(t)\int_\R e(\phi)dx+E_{(\phi,2,2)}(t)+\int_\R(|\psi|^2_\beta+|\nabla\psi|^4_\beta)dx\big),
\end{align}
where the positive constant \(C\) depends on \(N\) and the initial data.
\end{Prop}

\begin{proof}
Using the extrinsic version for Dirac-wave maps \eqref{phi-dwmap-extrinsic} we calculate
\begin{align*}
\frac{d}{dt}E_{(\phi,2,2)}(t)=&\int_\R\langle\partial_x(\Box\phi),\partial^2_{xt}\phi\rangle dx\\
=&\int_\R\langle\partial_x(\sff(d\phi,d\phi)),\partial^2_{xt}\phi\rangle dx 
+\int_\R\langle\partial_x(P(\sff(\psi,d\phi(\partial_t)),i\partial_t\cdot\psi)),\partial^2_{xt}\phi\rangle dx \\
&-\int_\R\langle\partial_x(P(\sff(\psi,d\phi(\partial_x)),i\partial_x\cdot\psi)),\partial^2_{xt}\phi\rangle dx.
\end{align*}
It is well-known that the contribution containing the second fundamental form can be estimated as
\[
|\langle\partial_x(\sff(d\phi,d\phi)),\partial^2_{xt}\phi\rangle|\leq C|d\phi|^3|d^2\phi|
\]
such that
\begin{align*}
\big\||d\phi|^3|d^2\phi|\big\|_{L^1(\R)}\leq & \big\||d\phi|^3\big\|_{L^2(\R)}\big\||d^2\phi|\big\|_{L^2(\R)} \\
\leq &C\big\|d|d\phi|^3\big\|_{L^\frac{2}{3}(\R)}\big\||d^2\phi|\big\|_{L^2(\R)} \\
\leq &C\big\||d\phi|^2|d^2\phi|\big\|_{L^\frac{2}{3}(\R)}\big\||d^2\phi|\big\|_{L^2(\R)} \\
\leq &C\big\||d\phi|\big\|^2_{L^2(\R)}\big\||d^2\phi|\big\|^2_{L^2(\R)},
\end{align*}
where we used the Sobolev embedding theorem in one dimension in the second step.
Regarding the other two terms we note
\begin{align*}
|\langle\partial_x(P(\sff(\psi,d\phi(\partial_t)),i\partial_t\cdot\psi)),\partial^2_{xt}\phi\rangle|
\leq & C(|\psi|_\beta^2|d\phi|^2|d^2\phi|+|\nabla\psi|_\beta|\psi|_\beta|d\phi||d^2\phi|+|\psi|_\beta^2|d^2\phi|^2) \\
\leq & C(|d^2\phi|^2+|\psi|_\beta^4|d\phi|^4+|\nabla\psi|^2_\beta|\psi|^2_\beta|d\phi|^2)\\
\leq & C(|d^2\phi|^2+|\psi|_\beta^2+|\nabla\psi|^4_\beta)
\end{align*}
and the remaining term can be estimated by the same method.
By adding up the various contributions we obtain the claim.
\end{proof}

\begin{Cor}
Let \((\phi,\psi)\colon\R^{1,1}\to\R^q\) be a Dirac-wave map.
Then the following estimate holds
\begin{align}
E_{(\phi,2,2)}(t)\leq Cf(t),
\end{align}
where the positive constant \(C\) depends on \(N\) and the initial data
and the function \(f(t)\) is finite for finite values of \(t\).
\end{Cor}
\begin{proof}
By the last Proposition we obtain an inequality of the form
\begin{align*}
\frac{d}{dt}E_{(\phi,2,2)}(t)\leq C(e^{Ct}E_{(\phi,2,2)}(t)+E_{(\phi,2,2)}(t)+e^{Ct}).
\end{align*}
The claim follows by the Gronwall-Lemma.
\end{proof}

At this point we have obtained sufficient control over the pair \((\phi,\psi)\) and its derivatives,
which is what we need to obtain a global solution of the system \eqref{phi-dwmap}, \eqref{psi-dwmap}.
In order to obtain a unique solution we will need the following

\begin{Prop}[Uniqueness]
\label{dwmap-prop-uniqueness}
Suppose we have two solutions of the Dirac-wave map system \eqref{phi-dwmap}, \eqref{psi-dwmap}.
If their initial data coincides, then they coincide for all times.
\end{Prop}

\begin{proof}
We make use of the extrinsic version of the Dirac-wave map system \eqref{phi-dwmap-extrinsic}, \eqref{psi-dwmap-extrinsic}.
Suppose that \(u,v\colon\R^{1,1}\to\R^q\) both solve \eqref{phi-dwmap-extrinsic} and \(\psi,\xi\in\Gamma(\Sigma\R^{1,1}\otimes\R^q)\)
both solve \eqref{psi-dwmap-extrinsic}, where \(\psi\) is a vector spinor along \(u\) and \(\xi\) a vector spinor along \(v\).
Since we have gained control over the \(L^2\)-norm of the second derivatives of \(\phi\), we have a pointwise bound on its
first derivatives. We set \(w:=u-v,~ \eta:=\psi-\xi\) and calculate
\begin{align*}
\frac{d}{dt}\frac{1}{2}\int_\R|dw|^2dx=&\int_\R\langle\frac{\partial w}{\partial t},\Box w\rangle dx \\
=&\int_\R\langle\frac{\partial w}{\partial t},\sff(u)(du,du)-\sff(v)(dv,dv)\rangle dx \\
&+\int_\R\epsilon_j\langle\frac{\partial w}{\partial t},P(u)(\sff(u)(e_j\cdot\psi,du(e_j)),i\psi)-P(v)(\sff(v)(e_j\cdot\xi,dv(e_j)),i\xi)\rangle dx.
\end{align*}
We rewrite the first term on the right hand side as follows
\[
\sff(u)(du,du)-\sff(v)(dv,dv)=\sff(u)(du,du)-\sff(v)(du,du)+\sff(v)(du+dv,dw).
\]
This yields
\[
\langle\frac{\partial w}{\partial t},(\sff(u)-\sff(v))(du,du)\rangle\leq C|du|^2|w||dw|.
\]
Using the orthogonality \(\partial_t\perp\sff\) we find
\[
\langle\frac{\partial w}{\partial t},\sff(v)(du,dw)\rangle=\langle\frac{\partial u}{\partial t},\sff(v)(du,dw)\rangle=\langle\frac{\partial u}{\partial t},(\sff(v)-\sff(u))(du,dw)\rangle\leq C|du|^2|w||dw|.
\]
The same argument also holds for the term involving \(dv\).
We rewrite 
\begin{align*}
P&(u)(\sff(u)(e_j\cdot\psi,du(e_j)),i\psi)-P(v)(\sff(v)(e_j\cdot\xi,dv(e_j)),i\xi) \\
=&(P(u)-P(v))(\sff(u)(e_j\cdot\psi,du(e_j)),i\psi)
+P(v)((\sff(u)-\sff(v))(e_j\cdot\psi,du(e_j)),i\psi)\\
&+P(v)(\sff(v)(e_j\cdot\eta,du(e_j)),i\psi)
+P(v)(\sff(v)(e_j\cdot\xi,dw(e_j)),i\psi)
+P(v)(\sff(v)(e_j\cdot\psi,du(e_j)),i\eta)
\end{align*}
such that we may estimate
\begin{align*}
|\langle\frac{\partial w}{\partial t},&P(u)(\sff(u)(e_j\cdot\psi,du(e_j)),i\psi)-P(v)(\sff(v)(e_j\cdot\xi,dv(e_j)),i\xi)\rangle| \\
&\leq C\big(|w||dw||\psi|_\beta^2|du|+|\eta|_\beta|w||du||\psi|_\beta+|w||dv||\eta|_\beta|\psi|_\beta\big)\leq C(|\eta|_\beta^2+|w|^2+|dw|^2).
\end{align*}
Regarding the spinors \(\psi,\xi\) we calculate
\begin{align*}
\frac{d}{dt}\frac{1}{2}\int_\R\langle\partial_t\cdot\eta,\eta\rangle dx=
\int_\R\epsilon_j\langle\sff(u)(e_j\cdot\psi,du(e_j))-\sff(v)(e_j\cdot\xi,dv(e_j)),\eta\rangle dx.
\end{align*}
We may rewrite the right hand side as follows
\begin{align*}
\sff(u)(e_j\cdot\psi,du(e_j))-\sff(v)(e_j\cdot\xi,dv(e_j))=&(\sff(u)-\sff(v))(e_j\cdot\psi,du(e_j))+\sff(v)(e_j\cdot\eta,dv(e_j)) \\
&+\sff(v)(e_j\cdot\xi,dw(e_j)).
\end{align*}
This allows us to obtain the following inequality
\begin{align*}
\frac{d}{dt}\int_\R|\eta|_\beta^2dx\leq C\int_\R(|\eta|_\beta|w|+|\eta|_\beta^2|dv|+|\eta|_\beta|dw||\xi|)dx
\leq C\int_\R(|\eta|_\beta^2+|w|^2+|dw|^2)dx.
\end{align*}
In total, we find
\begin{align*}
\frac{d}{dt}\int_\R(|\eta|_\beta^2+|w|^2+|dw|^2)dx\leq C\int_\R(|\eta|_\beta^2+|w|^2+|dw|^2)dx
\end{align*}
such that
\begin{align*}
\int_\R(|\eta|_\beta^2+|w|^2+|dw|^2)dx\leq e^{Ct}\int_\R(|\eta_0|_\beta^2+|w_0|^2+|dw_0|^2)dx.
\end{align*}
Consequently, we have \(u=v\) and \(\psi=\xi\) for all \(t\), whenever the initial data coincides.
\end{proof}

In the end we obtain the following existence result
\begin{Satz}
\label{satz-dwmap}
Let \(\R^{1,1}\) be two-dimensional Minkowski space and \((N,h)\) be a compact Riemannian manifold.
Then for any given initial data of the regularity
\begin{align*}
\phi(0,x)=&\phi_0(x)\in H^2(\R,N), \\
\frac{\partial\phi}{\partial t}(0,x)=&\phi_1(x)\in H^1(\R,TN),\\
\psi(0,x)=&\psi_0(x)\in H^1(\R,\Sigma\R^{1,1}\otimes\phi^\ast TN)\cap W^{1,4}(\R,\Sigma\R^{1,1}\otimes\phi^\ast TN)
\end{align*}
the equations for Dirac-wave maps \eqref{phi-dwmap}, \eqref{psi-dwmap} admit a global weak solution of the class
\begin{align*}
\phi\in H^2(\R^{1,1},N),\qquad \psi\in H^1(\R^{1,1},\Sigma\R^{1,1}\otimes\phi^\ast TN)\cap W^{1,4}(\R,\Sigma\R^{1,1}\otimes\phi^\ast TN),
\end{align*}
which is uniquely determined by the initial data.
\end{Satz}
\begin{proof}
By the energy estimates derived throughout this section we have enough control on the right hand sides of
\eqref{phi-dwmap-extrinsic} and \eqref{psi-dwmap-extrinsic} to extend the solution of the 
Dirac-wave map system for all times. The uniqueness follows from Proposition \ref{dwmap-prop-uniqueness}.
\end{proof}

\begin{Bem}
It is very desirable to get rid of the requirement 
\(\psi_0(x)\in  W^{1,4}(\R,\Sigma\R^{1,1}\otimes\phi^\ast TN)\) in Theorem \ref{satz-dwmap}.
However, in order to control the energy \(E_{(\phi,2,2)}(t)\) it seems to be necessary
to control \(L^p\) norms of \(\psi\) with \(p\geq 2\).
\end{Bem}

\begin{Bem}
The statement of Theorem \ref{satz-dwmap} still holds if we take into account an additional two-form
as was done in \cite{MR3305429} for the case of the domain being a closed Riemannian surface.
\end{Bem}

\par\medskip
\textbf{Acknowledgements:}
The author gratefully acknowledges the support of the Austrian Science Fund (FWF) 
through the START-Project Y963-N35 of Michael Eichmair
and the project P30749-N35 “Geometric variational problems from string theory”.
\bibliographystyle{plain}
\bibliography{mybib}

\end{document}